\documentclass[11pt]{article}
\usepackage[margin=1in]{geometry} 
\usepackage{amssymb,amsfonts,amsmath,bbm,mathrsfs,stmaryrd,mathtools,mathabx}
\usepackage{xcolor}
\usepackage{url}
\usepackage{relsize}
\usepackage{tikz-cd}
\usepackage[shortlabels]{enumitem}
\usepackage{braket}
\usepackage[T1]{fontenc}

\newcommand\xrsquigarrow[1]{%
  \mathrel{%
    \begin{tikzpicture}[baseline= {( $ (current bounding box.south) + (0,-0.5ex) $ )}]
      \node[inner sep=.5ex] (a) {$\scriptstyle #1$};
      \path[draw,
      <-,
      >=stealth,
      decorate,
      decoration={zigzag,amplitude=0.7pt,segment length=1.2mm,pre=lineto,pre length=4pt}] 
      (a.south east) -- (a.south west);
    \end{tikzpicture}}%
}

\usepackage[colorlinks,
             linkcolor=black!75!red,
             citecolor=blue,
             pdftitle={BU},
             pdfauthor={x},
             pdfproducer={pdfLaTeX},
             pdfpagemode=None,
             bookmarksopen=true
             bookmarksnumbered=true]{hyperref}

\usepackage{tikz}
\usetikzlibrary{arrows,calc,decorations.pathreplacing,decorations.markings,intersections,shapes.geometric,through,fit,shapes.symbols,positioning,decorations.pathmorphing}

\usepackage[amsmath,thmmarks,hyperref]{ntheorem}

\usepackage{comment}
\usepackage{cleveref}

\newcommand\nthalias[1]{\AddToHook{env/#1/begin}{\crefalias{lemma}{#1}}}

\nthalias{definition}
\nthalias{example}
\nthalias{examples}
\nthalias{remark}
\nthalias{remarks}
\nthalias{convention}
\nthalias{notation}
\nthalias{construction}
\nthalias{theoremN}
\nthalias{propositionN}
\nthalias{corollaryN}
\nthalias{lemma}
\nthalias{proposition}
\nthalias{corollary}
\nthalias{theorem}
\nthalias{conjecture}
\nthalias{question}
\nthalias{assumption}

\creflabelformat{enumi}{#2#1#3}

\crefname{section}{Section}{Sections}
\crefformat{section}{#2Section~#1#3} 
\Crefformat{section}{#2Section~#1#3} 

\crefname{subsection}{\S}{\S\S}
\AtBeginDocument{%
  \crefformat{subsection}{#2\S#1#3}%
  \Crefformat{subsection}{#2\S#1#3}% 
}

\crefname{subsubsection}{\S}{\S\S}
\AtBeginDocument{%
  \crefformat{subsubsection}{#2\S#1#3}%
  \Crefformat{subsubsection}{#2\S#1#3}% 
}

\theoremstyle{plain}

\newtheorem{lemma}{Lemma}[section]
\newtheorem{proposition}[lemma]{Proposition}
\newtheorem{corollary}[lemma]{Corollary}
\newtheorem{theorem}[lemma]{Theorem}

\theoremstyle{nonumberplain}

\theoremstyle{plain}
\theorembodyfont{\upshape}
\theoremsymbol{\ensuremath{\blacklozenge}}

\newtheorem{definition}[lemma]{Definition}
\newtheorem{example}[lemma]{Example}
\newtheorem{examples}[lemma]{Examples}
\newtheorem{remark}[lemma]{Remark}

\newtheorem{notation}[lemma]{Notation}

\crefname{definition}{definition}{definitions}
\crefformat{definition}{#2definition~#1#3} 
\Crefformat{definition}{#2Definition~#1#3} 

\crefname{ex}{example}{examples}
\crefformat{example}{#2example~#1#3} 
\Crefformat{example}{#2Example~#1#3} 

\crefname{exs}{example}{examples}
\crefformat{examples}{#2example~#1#3} 
\Crefformat{examples}{#2Example~#1#3} 

\crefname{remark}{remark}{remarks}
\crefformat{remark}{#2remark~#1#3} 
\Crefformat{remark}{#2Remark~#1#3} 

\crefname{remarks}{remark}{remarks}
\crefformat{remarks}{#2remark~#1#3} 
\Crefformat{remarks}{#2Remark~#1#3} 

\crefname{convention}{convention}{conventions}
\crefformat{convention}{#2convention~#1#3} 
\Crefformat{convention}{#2Convention~#1#3} 

\crefname{notation}{notation}{notations}
\crefformat{notation}{#2notation~#1#3} 
\Crefformat{notation}{#2Notation~#1#3} 

\crefname{claim}{claim}{claims}
\crefformat{claim}{#2claim~#1#3} 
\Crefformat{claim}{#2Claim~#1#3} 

\crefname{conjecture}{conjecture}{conjectures}
\crefformat{conjecture}{#2conjecture~#1#3} 
\Crefformat{conjecture}{#2Conjecture~#1#3}

\crefname{lemma}{lemma}{lemmas}
\crefformat{lemma}{#2lemma~#1#3} 
\Crefformat{lemma}{#2Lemma~#1#3} 

\crefname{proposition}{proposition}{propositions}
\crefformat{proposition}{#2proposition~#1#3} 
\Crefformat{proposition}{#2Proposition~#1#3} 

\crefname{question}{question}{questions}
\crefformat{question}{#2question~#1#3} 
\Crefformat{question}{#2Question~#1#3}

\crefname{corollary}{corollary}{corollaries}
\crefformat{corollary}{#2corollary~#1#3} 
\Crefformat{corollary}{#2Corollary~#1#3} 

\crefname{theorem}{theorem}{theorems}
\crefformat{theorem}{#2theorem~#1#3} 
\Crefformat{theorem}{#2Theorem~#1#3} 

\crefname{enumi}{}{}
\crefformat{enumi}{#2#1#3}
\Crefformat{enumi}{#2#1#3}

\crefname{assumption}{assumption}{Assumptions}
\crefformat{assumption}{#2assumption~#1#3} 
\Crefformat{assumption}{#2Assumption~#1#3} 

\crefname{equation}{}{}
\crefformat{equation}{(#2#1#3)} 
\Crefformat{equation}{(#2#1#3)}

\numberwithin{equation}{section}

\theoremstyle{nonumberplain}
\newtheorem{proof}{Proof}

\newcommand\bC{{\mathbb C}}

\newcommand\bN{{\mathbb N}}

\newcommand\bR{{\mathbb R}}

\newcommand\bZ{{\mathbb Z}}

\newcommand\cC{{\mathcal C}}

\newcommand\cQ{{\mathcal Q}}

\newcommand\con{{\rm cone}}

\newcommand\restr[2]{{% we make the whole thing an ordinary symbol
  \left.\kern-\nulldelimiterspace % automatically resize the bar with \right
  #1 % the function
  \vphantom{\big|} % pretend it's a little taller at normal size
  \right|_{#2} % this is the delimiter
  }}

\DeclareMathOperator{\id}{id}

\newcommand{\cat}[1]{\textsc{#1}}

\newcommand{\qedhere}{\mbox{}\hfill\ensuremath{\blacksquare}}

\def\polhk#1{\ssetbox0=\hbox{#1}{\ooalign{\hidewidth
    \lower1.5ex\hbox{`}\hidewidth\crcr\unhbox0}}}

\newcommand{\bes}{\begin{equation*}}
\newcommand{\ees}{\end{equation*}}
\newcommand{\be}{\begin{equation}}
\newcommand{\ee}{\end{equation}}

\begin{document}
\date{}

\newcommand{\Addresses}{{% additional braces for segregating \footnotesize
  \bigskip
  \footnotesize

  \textsc{Department of Mathematics, University at Buffalo}
  \par\nopagebreak
  \textsc{Buffalo, NY 14260-2900, USA}  
  \par\nopagebreak
  \textit{E-mail address}: \texttt{achirvas@buffalo.edu}

  \medskip

  \textsc{Instytut matematyczny, Uniwersytet Wroc{\l}awski}
  \par\nopagebreak
  \textsc{PL. Grunwaldzki 2, Wroc{\l}aw, 50-384 Poland}  
  \par\nopagebreak
  \textit{E-mail address}: \texttt{mariusz.tobolski@math.uni.wroc.pl}

}}

\title{Non-commutative classifying spaces of groups\\
via quasi-topologies and pro-$C^*$-algebras}
\author{Alexandru Chirvasitu and Mariusz Tobolski}

\maketitle
\begin{abstract}
  For a locally compact Hausdorff topological group $G$, we construct a universal pro-$C^*$-algebra $C(E^+G)$ as the non-commutative geometer's analogue of the total space $EG$ of the classifying principal $G$-bundle $EG\to BG$. The pro-$C^*$-algebra $C(EG)$ of (possibly unbounded) continuous functions on $EG$ is then recoverable as the abelianization of $C(E^+G)$. Along the way, we develop various aspects of the theory of quasi-topological $G$-spaces and $G$-pro-$C^*$-algebras.
\end{abstract}


\vspace{.5cm}


%\tableofcontents

%%%%%%%%%%%%%%%%%%%%%%%%%%%%%%%%%%%%%%%%%%%%%%%%%%%%%%%%%%%%%%%%%%%%%%%%%%%%%%%%%%%%%%
%%%%%%%%%%%%%%%%%%%%%%%%%%%%%%%%%%%%%%%%%%%%%%%%%%%%%%%%%%%%%%%%%%%%%%%%%%%%%%%%%%%%%%

\section*{Introduction}

Principal bundles are the currency of gauge theory, and as such are ubiquitous in the modern mathematical-physics literature (see e.g.~\cite{isham,naber}). Recall \cite[\S 14.4]{t-td08} that for a topological group $G$ one defines the {\em classifying space} $BG$ as the base space of the {\em universal bundle} $EG\to BG$, universal in the following sense: for any other {\it numerable} \cite[\S 14.3]{t-td08} principal $G$-bundle $X\to M$ there exists a $G$-map $X\to EG$ which is unique up to $G$-homotopy.

The construction of $EG$ for an arbitrary topological group is due to Milnor~\cite{j-m56}. $BG$ classifies all numerable principal $G$-bundles in the sense that there is a bijection between isomorphism classes of numerable principal $G$-bundles $X\to M$ and homotopy classes of maps $M\to BG$. The universal space $EG$ can be used to define equivariant (co)homology of a space \cite[\S 2]{ab_equiv}, while the theory of characteristic classes in many cases amounts to computing the cohomology groups of $BG$. The importance of $BG$ in physics is difficult to overstate (as the already-cited \cite{isham,naber} and their ample references will confirm) but we mention, for instance, that it features prominently in Dijkgraaf--Witten theory~\cite{dw90}.

Gauge theory is also central to Connes' {\it non-commutative geometry} programme ~\cite{connes}. Here, the starting point is {\it Gelfand--Naimark duality}~\cite{gn-43}, which provides a contravariant equivalence between the category of locally compact Hausdorff spaces and that of commutative $C^*$-algebras. For a locally compact Hausdorff group $G$ this equivalence can be extended to locally compact Hausdorff $G$-spaces and commutative $G$-$C^*$-algebras.

Many properties of group actions on locally compact Hausdorff spaces have non-commutative analogs, in the realm of group actions on $C^*$-algebras. The freeness of an action, for instance, can be generalized in many ways (e.g.~\cite{ncp-09} for actions of finite groups). An analog of a numerable principal $G$-bundle in this context was introduced by the second author in~\cite{mt22}, generalizing the notion of the local-triviality dimension~\cite{ghtw-18} (with a topological structure group).

There is some prior work by M\"uller~\cite{a-m92} and \DJ ur\dj evi\'c~\cite{m-d97} on non-commutative (analogues of) classifying spaces, taking different approaches in somewhat different settings than here: both authors work with actions of {\it compact matrix quantum groups} \cite{wor_cmpg} on unital $C^*$-algebras. 

M\"uller's elegant machinery is only applicable to continuous fields of $C^*$-algebras and compact quantum groups with characters.  Because the resulting classifying spaces are intimately linked to the {\it classical} compact groups of characters of the underlying compact quantum groups, only a ``classical shadow'' of the non-triviality of non-commutative principal bundle considered in his theory is visible. For Woronowicz's $SU_q(n)$~\cite{woronowicz}, for instance, the largest classical subgroup is (for $q\in (0,1)$ sufficiently close to $1$ \cite[Proposition 17]{a-m92}) the maximal torus of $SU(n)$ that ``survives'' the $q$-deformation; the passage from a compact quantum group to its maximal classical subgroup can thus be rather drastic.

\DJ ur\dj evi\'c's approach relies on the crucial observation familiar from algebraic topology that a free action of a compact Lie group $G$ on a compact Hausdorff $X$ space always gives rise to a principal $G$-bundle $X\to X/G$. Now, if $H$ is a closed subgroup of compact Hausdorff group $G$, then $EG$ is the universal space for $H$ if and only if $G\to G/H$ is a principal $H$-bundle. Since this is always the case for some embedding $H\le U(m)$ into a unitary group for some $m\geq 0$, the universal space of any compact Lie group coincides with some $EU(m)$.

In the non-commutative case there is no guarantee that free actions of compact matrix quantum groups give rise to locally trivial non-commutative principal bundles. In fact,~\cite{cpt-21} shows that while the usual $U(1)$-action on the Cuntz algebra~$\mathcal{O}_2$~\cite{cuntz} is free, it cannot define a locally trivial principal bundle in the sense of~\cite{ghtw-18} and~\cite{mt22}. Therefore, the work of \DJ ur\dj evi\'c does not generalize the notion of a locally trivial principal bundle, but a condition that is strictly stronger.

The classifying space of a group is typically not locally compact: $B\left(\bZ/2\bZ\right)$, for instance, is the infinite real projective space. For that reason, a non-commutative approach to the topic must go beyond the theory of $C^*$-algebras. Here, we take a cue from the fact that the category of {\it quasi-topological spaces} of Spanier~\cite{e-s63} is equivalent to the category of commutative unital {\it pro-$C^*$-algebras}~\cite{phil1} via the Gelfand--Naimark-type duality of Dubuc and Porta~\cite{dubuc}. Taking that equivalence as motivation on the one hand and validation on the other, we construct the non-commutative classifying space of a completely Hausdorff group equipped with a quasi-topology using universal pro-$C^*$-algebras~\cite{ta-l10,phil2}.

The construction generalizes a slightly altered Milnor space $EG$. Because the category of completely Hausdorff quasi-topological spaces is equivalent \cite[Proposition 2.6]{phil1} to the category of (again completely Hausdorff) spaces with distinguished families of compact subsets, our construction applies to all locally compact Hausdorff groups. We then show that the aforementioned construction gives a contractible pro-$C^*$-algebra in a suitable sense, and use the universal principal bundle thus obtained to generalize the definition of a~non-commutative numerable principal bundle given in~\cite{mt22}.

The paper is organized as follows.

\Cref{se:prel} collects basic definitions and results regarding quasi-topological spaces and pro-$C^*$-algebras.

In~\Cref{se:quact}, we extend the anti-equivalence of quasi-topological space and commutative unital pro-C*-algebras to the equivariant setting.

\Cref{se:nus} introduces the notion of a non-commutative classifying space of a completely Hausdorff quasi-topological group. First, we slightly change the Milnor construction: our $EG$ is different as a set while carrying the same topology. We show that we still obtain a universal bundle and generalize the thus obtained space using the language of universal pro-$C^*$-algebras. Along the way we reformulate the definition of the pro-$C^*$-relation to one better suited for our purposes.

In~\Cref{se:ncprin}, we use our construction to consider pro-$C^*$-algebraic and $C^*$-algebraic analogs of principal bundles, which generalize many notions known from the literature.

Although the results of this paper are true for {\it group objects} in the category of quasi-topological spaces (or groups {\it internal} to it) (see~\cite{catm-25}), we only deal with locally compact Hausdorff topological groups with their naturally induced quasi-topology.

%%%%%%%%%%%%%%%%%%%%%%%%%%%%%%%%%%%%%%%%%%%%%%%%%%%%%%%%%%%%%%%%%%%%%%%%%%%%%
\subsection*{Acknowledgments}

A.C. was partially supported by NSF grant DMS-2001128. M.T. was supported by NCN grant 2020/36/C/ST1/00082 entitled {\it Non-Commutative universal spaces for groups and quantum groups}. This work was completed when M.T. was visiting University at Buffalo, and he is very grateful to the University for excellent working conditions and amazing hospitality.

We are also grateful for the anonymous referee's valuable suggestions towards improving initial draft versions of the paper. 

%%%%%%%%%%%%%%%%%%%%%%%%%%%%%%%%%%%%%%%%%%%%%%%%%%%%%%%%%%%%%%%%%%%%%%%%%%%%%%%%%%%%%%
%%%%%%%%%%%%%%%%%%%%%%%%%%%%%%%%%%%%%%%%%%%%%%%%%%%%%%%%%%%%%%%%%%%%%%%%%%%%%%%%%%%%%%
\section{Preliminaries}\label{se:prel}

Here we collect some basic definitions and results on quasi-topologies and pro-$C^*$-algebras.

%%%%%%%%%%%%%%%%%%%%%%%%%%%%%%%%%%%%%%%%%%%%%%%%%%%%%%%%%%%%%%%%%%%%%%%%%%%%%

\subsection{Quasi-topological spaces}

First, we recall the definition of a quasi-topological space due to Spanier~\cite{e-s63}.

\begin{definition}\label{def:qtop}
  A {\em quasi-topological space} (occasionally {\it quasi-space} for short) is a pair $(X,Q)$, where $X$ is a set and $Q$ is an assignment, called {\em quasi-topology}, which to each compact Hausdorff space $K$ assigns a set $Q(K,X)$ of maps from $K$ to $X$ such that
  \begin{enumerate}[{(Q}1{)}]
  \item\label{item:q1} All constant maps from $K$ to $X$ are in $Q(K,X)$.
  \item\label{item:q2} If $\varphi:K_1\to K_2$ is a continuous map between compact Hausdorff spaces, and $\psi\in Q(K_2\,,X)$, then $\psi\circ \varphi\in Q(K_1,X)$. 
  \item\label{item:q3} If $\varphi:K_1\to K_2$ is a continuous surjective map between compact Hausdorff spaces, and $\psi:K_2\to X$, then $\psi\in Q(K_2,X)$ if and only if $\psi\circ \varphi\in Q(K_1\,,X)$.
  \item\label{item:q4} If $K$ is a disjoint union of compact Hausdorff spaces $K_1$ and $K_2$, then $\varphi\in Q(K,X)$ if and only if $\varphi|_{K_i}\in Q(K_i\,,X)$, $i=1,2$.
  \end{enumerate}
  A map $f:X\to Y$ between quasi-topological spaces is called {\em quasi-continuous} if for every $\varphi\in Q(K,X)$ we have that $f\circ \varphi\in Q(K,Y)$. A {\em quasi-homeomorphism} is a bijective quasi-continuous map whose inverse is quasi-continuous.
\end{definition}

\begin{remark}\label{re:daysays}
  Condition \ref{item:q4} of \Cref{def:qtop} appears on \cite[p.3]{bj-d68} as Q3, and is phrased somewhat differently: Day works with arbitrary finite unions, including the empty union. The latter implies in particular that the unique map $\emptyset\to X$ belongs to $Q(\emptyset,X)$.

  That \ref{item:q4} recovers the same content is easy to see: unions of more than two sets come about inductively, whereas the fact that $\emptyset\to X$ is distinguished follows, for instance, by taking $K_1$ to be a singleton and $K_2$ empty in \ref{item:q4}.
\end{remark}

There are several induced quasi-topologies~(see \cite[Section~3]{e-s63}):

\begin{examples}\label{exs:quasitops}
  \begin{enumerate}[(1)]
 
  \item\label{item:subspace} Let $Y$ be a subset of a quasi-topological space $(X,Q)$. Then the {\em subspace quasi-topology} $Q_{\subseteq}$ on $Y$ is defined by the condition that for every compact Hausdorff space $K$ a map $\varphi\in Q_\subseteq(K,Y)$ if and only if the composition of $\varphi$ and the inclusion $Y\subseteq X$ belongs to $Q(K,X)$.
  \item Let $\{(X_\alpha,Q_\alpha)\}$ be a family of quasi-topological spaces indexed by a set $I$. Then the {\em product quasi-topology} $Q_\times$ on $X:=\prod X_\alpha$ is defined by the condition that for every compact Hausdorff space $K$ a map $\varphi\in Q_\times(K,X)$ if and only if $\pi_\alpha\circ \varphi\in Q_\alpha(K,X_\alpha)$, where $\pi_\alpha:X\to X_\alpha$ is the canonical projection map.
  \item\label{item:quot.qtop} Let $\pi:X\to Y$ be a surjection from a quasi-topological space $(X,Q)$ to a set. Then the {\em quotient quasi-topology} $Q_\pi$ on $Y$ is defined by the condition that for every compact Hausdorff space $K$ a~map $\varphi\in Q_\pi(K,Y)$ if and only if there exists a compact Hausdorff space $K'$, a~continuous surjective map $\pi':K'\to K$, and $\varphi'\in Q(K',X)$ rendering the diagram
    \begin{equation*}
      \begin{tikzcd}
        K' \arrow{d}[swap]{\pi'} \arrow{r}{\varphi'} & X \arrow{d}{\pi}\\
        K \arrow{r}{\varphi} & Y
      \end{tikzcd}
    \end{equation*}
    commutative.
  \end{enumerate}
\end{examples}

We will make use of the following two lemmas proved by Spanier~\cite[Lemma~3.1, Lemma~3.2]{e-s63}.
\begin{lemma}\label{quotlem}
Let $\pi:X\to Y$ be a surjection from a quasi-topological space onto a set and let $Y$ be given a quotient quasi-topology. Then a map $f:Y\to Z$ from $Y$ to a quasi-topological space $Z$ is quasi-continuous if and only if $f\circ \pi:X\to Z$ is quasi-continuous.
\end{lemma}
\begin{lemma}\label{prodquot}
Let $\{\pi_\alpha\}$ be a family of surjective maps from quasi-topological spaces $\{X_\alpha\}$ to sets $\{Y_\alpha\}$. If $\pi:X\to Y$ is the surjection from the product of the quasi-topological spaces $\{X_\alpha\}$ to the product of sets $\{Y_\alpha\}$ given by the family $\{\pi_\alpha\}$, then the product quasi-topology on $Y$ coincides with the quotient quasi-topology given by $\pi$.
\end{lemma}

Every topological space $X$ can be turned into a quasi-topological space by taking for $Q(K,X)$ the set $C(K,X)$ of {\it all} continuous maps from $K$ to $X$, but it is worth noticing that not every quasi-topology arises that way~(see~\cite[Lemma~5.5]{e-s63}). In this context, let us recall a definition due to Phillips~\cite[Definition~2.5, Proposition~2.6]{phil1}.

\begin{definition}\label{def:compfam}
  A {\em distinguished family of compact subsets} $\mathscr{F}_X$ on a topological space $X$ is a family of compact subsets of $X$ such that
  \begin{enumerate}[{(K}1{)}]
  \item\label{item:k1} $\{x\}\in\mathscr{F}_X$ for all $x\in X$.
  \item\label{item:k2} If $F\in\mathscr{F}_X$ and $K$ is a compact subset of $F$, then $K\in\mathscr{F}_X$.
  \item\label{item:k3} If $F_1$, $F_2\in\mathscr{F}_X$, then $F_1\cup F_2\in\mathscr{F}_X$.
  \item\label{item:k4} A subset $Y$ of $X$ is closed if and only if $Y\cap F$ is closed for all $F\in\mathscr{F}_X$.
  \end{enumerate}
  A {\em morphism} $f:(X,\mathscr{F}_X)\to (Y,\mathscr{F}_Y)$ between spaces with distinguished families of compact subsets is a continuous map $f:X\to Y$ such that $f(F)\in\mathscr{F}_Y$ for all $F\in\mathscr{F}_X$.
\end{definition}
Note that compactly generated spaces (in particular, locally compact spaces) are spaces with distinguished families of compact subsets. The following class of (quasi-)topological spaces is important from the point of view of pro-$C^*$-algebra theory (see~\cite[discussion before~Theorem~4.11]{dubuc} and~\cite[Definition~2.2]{phil1}).
\begin{definition}\label{def:comphaus}
  A (quasi-)topological space $X$ is {\it completely Hausdorff} if for any two different points $x, y\in X$ there is a (quasi-)continuous function $f:X\to[0,1]$ such that $f(x)=0$ and \mbox{$f(y)=1$}. Equivalently, $X$ is completely Hausdorff if for any two different points $x, y\in X$ there is a \mbox{(quasi-)}continuous function $f:X\to\bR$ such that $f(x)\neq f(y)$.
\end{definition}
 
%%%%%%%%%%%%%%%%%%%%%%%%%%%%%%%%%%%%%%%%%%%%%%%%%%%%%%%%%%%%%%%%%%%%%%%%%%%%%

\subsection{Pro-$C^*$-algebras}

Pro-$C^*$-algebras are also known under the names LM$C^*$-algebras and b*-algebras. We refer the reader to the book of Fragoulopoulou~\cite[p.~99--101]{m-f05} for a historical review of the theory which goes back to 1952. \cite{phil1,phil2} are good sources for the material in the present subsection.

\begin{definition}\label{def:proca}
A {\em pro-$C^*$-algebra} $A$ is a complete Hausdorff topological $*$-algebra over $\mathbb{C}$ whose topology is determined by a family of continuous $C^*$-seminorms $\{p_\lambda\}_{\lambda\in\Lambda}$. Morphisms between pro-$C^*$-algebras are continuous $*$-homomorphisms. If the family of $C^*$-seminorms defining the topology on a pro-$C^*$-algebra is countable, we say that it is a {\em$\sigma$-$C^*$-algebra}.
\end{definition}
We denote the set of all continuous $C^*$-seminorms on a pro-$C^*$-algebra $A$ by $S(A)$ and write $\{p_\lambda\}_{\lambda\in\Lambda}\subseteq S(A)$. In the case of $\sigma$-$C^*$-algebras we can always chose $\Lambda=\bN$ in~\Cref{def:proca} and make $*$-homomorphisms $A\to A_\lambda$ surjective. Note that in contrast to the theory of $C^*$-algebras, $*$-homo\-morphisms of pro-$C^*$-algebras might not be continuous. However, $*$-homomorphisms between $\sigma$-$C^*$-algebras are automatically continuous~\cite[Theorem~5.2]{phil1}. Any $C^*$-algebra is a pro-$C^*$-algebra in the above sense. Many features of the theory of \mbox{$C^*$-alge}\-bras (functional calculus, tensor products, multiplier algebras, Hilbert modules, etc.) can be easily extended to the theory of pro-$C^*$-algebras due to the following result (\cite[Proposition 1.2]{phil1}).
\begin{proposition}
Every pro-$C^*$-algebra is an inverse limit of $C^*$-algebras in the category of topological $*$-algebras and continuous $*$-homomorphisms, namely $A\cong\varprojlim A_\lambda$, where $A_\lambda:=A/\ker p_\lambda$.
\end{proposition}

%%%%%%%%%%%%%%%%%%%%%%%%%%%%%%%%%%%%%%%%%%%%%%%%%%%%%%%%%%%%%%%%%%%%%%%%%%%%%

\subsection{Categorical (anti-)equivalences}\label{subse:catres}

In general, for a category $\cC$, we write $\cC(c,d)$ for the set of morphisms between two objects $c,d$ thereof. We will the following categories:
 \begin{notation}\label{not:cats}
  We will often use the following categories:
  \begin{enumerate}[(1)]
  \item $\cat{Top}$ is the category of topological spaces, and $\cat{Top}_{cH}\subset \cat{Top}$ the  subcategory of completely Hausdorff (see~\Cref{def:comphaus}) ones.
  \item $\cat{Top}_{\mathscr{F}}$ is the category of topological spaces with distinguished families of compact subsets (see~\Cref{def:compfam}), and $\cat{Top}_{cH,\mathscr{F}}\subset \cat{Top}_{\mathscr{F}}$ the subcategory of completely Hausdorff ones. 
   \item $\cQ$ is the category of quasi-topological spaces (see~\Cref{def:qtop}) and $\cQ_{cH}\subset \cQ$ the subcategory of completely Hausdorff ones.  
  \item\label{item:allcast} The generic symbol $\cC^*$ denotes one of the several categories of (pro-)$C^*$-algebras, and we distinguish between these by additional subscripts: a `c' for commutative, a `1' for unital, a `pro' for pro-$C^*$-algebras that may or may not carry the additional structure. Thus:
    \begin{itemize}
    \item $\cC^*$ alone means plain $C^*$-algebras, possibly non-unital, with continuous morphisms;
    \item $\cC^*_{pro,1}$ is the category of unital pro-$C^*$-algebras;
    \item $\cC^*_{pro,c}$ that of possibly non-unital commutative pro-$C^*$-algebras;
    \item $\cC^*_{c,1}$ that of commutative unital $C^*$-algebras, etc. etc.
    \end{itemize}
    We will occasionally phrase a statement uniformly for several of the categories. We indicate this by a blanket stand-in symbol:  $\cC^*_{\bullet}$ means any of the above or $\cC^*_{pro,\bullet}$ means pro-$C^*$-algebras, plain, commutative, unital or both.
    \end{enumerate}   
\end{notation}

\begin{remark}\label{rem:q}
  The symbol `$\cQ$' is meant as reminiscent of Day's notation \cite[pp.2-3]{bj-d68} for the category of quasi-topological spaces. We might occasionally trail that source's notation in other ways.
\end{remark}

We review various equivalences between the aforementioned categories. The following result is due to Phillips~\cite[Proposition~2.6]{phil1}: the cited result's somewhat different phrasing notwithstanding, its proof does deliver the apparently stronger statement.

\begin{proposition}\label{pr:chqt}
  There is an equivalence between the categories $\cat{Top}_{cH,\mathscr{F}}$ and $\cQ_{cH}$ via the functor which assigns to each space with a distinguished family of compact subsets $(X,\mathscr{F}_X)$ the quasi-topology given by the condition that for any compact Hausdorff space $K$ a map $\varphi\in Q(K,X)$ if and only if $\varphi$ is continuous and $\varphi(K)\in\mathscr{F}_X$.
\end{proposition}

% \Cref{pr:chqt} immediately implies:
% \begin{corollary}\label{cor:chqt}
%   The categories $\cat{Top}_{cH,\mathscr{F}}$ and $\cQ_{cH}$ are equivalent.
% \end{corollary}

The following result of Dubuc and Porta~\cite[Theorem~4.11(e)]{dubuc} (reproved by Phillips in~\cite[Theorem~2.7]{phil1}) plays the role of the Gelfand-Naimark theorem and allows us to treat unital pro-$C^*$-algebras as algebras of functions on {\em (completely Hausdorff) non-commutative quasi-topological spaces}.
\begin{theorem}\label{th:antieq}
  There is an anti-equivalence between the categories $\cQ_{cH}$ and $\cC^*_{pro,c,1}$ via the functor 
  \begin{equation*}
    X\longmapsto C(X),
  \end{equation*}
  where $C(X)$ denotes the $*$-algebra of quasi-continuous complex-valued functions with the topology given by $C^*$-seminorms $p_{K,\varphi}(f):=\sup_{x\in K}|f(\varphi(x))|$ for any compact Hausdorff space $K$ and $\varphi\in Q(K,X)$. The inverse functor is denoted by $A\mapsto {\rm Spec}(A)$.
\end{theorem}

%%%%%%%%%%%%%%%%%%%%%%%%%%%%%%%%%%%%%%%%%%%%%%%%%%%%%%%%%%%%%%%%%%%%%%%%%%%%%
%%%%%%%%%%%%%%%%%%%%%%%%%%%%%%%%%%%%%%%%%%%%%%%%%%%%%%%%%%%%%%%%%%%%%%%%%%%%%

%%%%%%%%%%%%%%%%%%%%%%%%%%%%%%%%%%%%%%%%%%%%%%%%%%%%%%%%%%%%%%%%%%%%%%%%%%%%%
\section{Group actions on commutative pro-$C^*$-algebras}\label{se:quact}

In this section, we extend the anti-equivalence in~\Cref{th:antieq} to the  setting of $G$-actions on quasi-topological spaces and commutative pro-C*-algebras.

A {\em quasi-continuous action} of a locally compact group $G$ on a~quasi-topological space $X$ is a~quasi-continous map $G\times X\to X$, where on the domain we consider the product quasi-topology, satisfying the usual axioms of an action. We also follow~\cite[Definition~5.1]{nc-p89} and say that $G$ {\em acts jointly continuously} on a pro-C*-algebra $A$ if there exists a group homomorphism $\alpha:G\to {\rm Aut}(A)$, where ${\rm Aut}(A)$ denotes the continuous $*$-automorphisms of $A$, such that the map
\[
G\times A\ni (g,a)\longmapsto \alpha_g(a):=(\alpha(g))(a)\in A
\]
is continuous for the two obvious topologies.

The following result gives the desired correspondence between actions of locally compact Hausdorff groups on commutative pro-$C^*$-algebras and their respective spectra.

\begin{theorem}\label{th:act}
  Let $G$ be a locally compact Hausdorff group, $X$ a quasi-topological space and $A:=C(X)$ its corresponding unital pro-$C^*$-algebra. For an action
  \begin{equation*}
    G\times A\xrightarrow{\quad\alpha\quad}A
  \end{equation*}
  of $G$ on $A$, the following conditions are equivalent:
  \begin{enumerate}[(a)]

  \item\label{item:space} $\alpha$ is induced by a quasi-continuous action $G\times X\xrightarrow{\quad}X$;

  \item\label{item:joint} $\alpha$ is a jointly continuous action of $G$ on $A$;

  \item\label{item:strong} point-uniform continuity, in the sense that for every $a\in A$ the map
    \begin{equation*}
      G\ni g\xmapsto{\quad}ga\in A
    \end{equation*}
    is continuous.
  
  \end{enumerate}
\end{theorem}
\begin{proof}
  We prove three implications to complete the circular chain.
  \begin{enumerate}[]

  \item {\bf \Cref{item:space} $\Longrightarrow$ \Cref{item:joint}:} Recall from the proof of \cite[Theorem 2.7]{phil1} that $X$ is the space of continuous $*$-morphisms $A\to \bC$, and its distinguished family of compact sets \cite[Definition 2.5]{phil1} consists of the (compact) spaces $X_p\subseteq X$ of morphisms $A_p\to \bC$ where
    \begin{itemize}
    \item $p$ is a $C^*$ seminorm on $A$;
    \item and $A\to A_p$ is the corresponding $C^*$ quotient (the quotient by the kernel of $p$ being automatically complete \cite[Corollary 1.12]{phil1}).
    \end{itemize}
    The requirement that the action $G\times X\to X$ (which we typically depict by juxtaposition) be a morphism of quasi-topological spaces (or, equivalently \cite[Proposition 2.6]{phil1}, one of completely Hausdorff spaces with distinguished compact-set families) means simply that
    \begin{itemize}
    \item $(g,x)\xmapsto{} gx$ is continuous;
    \item and it maps sets of the form
      \begin{equation*}
        K\times X_p,\quad K\subseteq G\text{ compact}
      \end{equation*}
      into $X_q$ (varying $q$ for varying $p$, etc.).
    \end{itemize}
    The action on $A$ induced by that on $X$ is
    \begin{equation*}
      \alpha(g,f):=f(g^{-1}\cdot),\ \forall g\in G,\ f\in A=C(X).
    \end{equation*}
    Consider nets $g_{\lambda}\to 1$ in $G$ and $f_{\lambda}\to f$ in $A$ respectively. The latter convergence means uniform convergence on every distinguished compact $X_p$. Because $G$ is locally compact, $g_{\lambda}$ can be assumed inside a compact neighborhood $K$ of $1\in G$.

    Consider, now, a $C^*$ seminorm $p$ with its attached compact set $X_p\subseteq X$. We noted above that the $G$-action on $X$ (co)restricts to
    \begin{equation*}
      K\times X_p\xrightarrow{\quad}X_q\text{ for some }q,
    \end{equation*}
    and the uniform convergence $f_{\lambda}\to f$ on $X_q$ then implies the uniform convergence
    \begin{equation*}
      \alpha(g_{\lambda},f_{\lambda}) = f_{\lambda}(g_{\lambda}^{-1}\cdot) \to f(1\cdot) = f
    \end{equation*}
    on $X_p$. This finishes the proof of the present implication.

  \item {\bf \Cref{item:joint} $\Longrightarrow$ \Cref{item:strong}:} This is obvious, the second property clearly being formally weaker.

  \item {\bf \Cref{item:strong} $\Longrightarrow$ \Cref{item:space}:} Note first that for any automorphism $\theta$ of $A$ as a $*$-algebra and any $C^*$ seminorm $p$, the twist
    \begin{equation*}
      {}^{\theta}p:=p(\theta^{-1}\cdot)
    \end{equation*}
    is again a $C^*$ seminorm. This implies, in particular, that any group action on $A$ (regardless of any continuity requirements) induces one on the space $X$ of $*$-algebra morphisms $A\to \bC$ continuous with respect to the family of all $C^*$ seminorms. This allows us to speak of {\it the} action of $G$ on $X$ (induced by that on $A$), as we will in the sequel.

    Two claims are now implicit in the goal; we tackle them in turn.

    \begin{enumerate}[(I)]

    \item\label{item:kpq} {\bf The action $G\times X\to X$ respects families of distinguished compact subsets.} In other words, preserving the notation in the proof of the implication \Cref{item:space} $\Longrightarrow$ \Cref{item:joint}, we have to argue that the action factors as
      \begin{equation}\label{eq:kpq}
        K\times X_p\to X_q\text{ for some $C^*$ seminorm $q$}
      \end{equation}
      for every $C^*$ seminorm $p$ and compact $K\subseteq G$. Since $X_q$ is precisely the space of characters on $A$ continuous with respect to $q$, this amounts to showing that for every compact $K\subseteq G$ and $C^*$ seminorm $p$ on $A$ there is a $C^*$ seminorm $q$ dominating all of the twists
      \begin{equation*}
        {}^gp = p(g^{-1}\cdot),\ g\in K. 
      \end{equation*}
      To see this, simply set
      \begin{equation*}
        q:=\sup_{g\in K}{}^gp = \sup_{g\in K}p(g^{-1}\cdot).
      \end{equation*}
      The only possible issue is that of finiteness: is it the case that
      \begin{equation*}
        \sup_{g\in K}p(g^{-1}a)<\infty,\ \forall a\in A?
      \end{equation*}
      The affirmative answer follows from the assumed continuity property of the action on $A$: for every $a\in A$ the map
      \begin{equation*}
        K\ni g\xmapsto{\quad}g^{-1}a\in A
      \end{equation*}
      is assumed continuous, further composition with $p$ will provide a continuous map $K\to \bR_{\ge 0}$, and {\it that} map must be bounded because its domain is compact.

    \item {\bf The action $G\times X\to X$ is continuous.} Recall from the proof of \cite[Theorem 2.7]{phil1} that the topology on $X$ is in fact the direct-limit topology induced by
      \begin{equation*}
        X = \varinjlim_p X_p,\quad p\text{ ranging over $C^*$ seminorms on $A$}
      \end{equation*}
      (the connecting maps $X_p\to X_q$ being the inclusions resulting from domination relations $p\le q$).

      Now because $G$ is locally compact, the canonical morphism
      \begin{equation}\label{eq:xtimesvarinjlim}
        \varinjlim_p \left(G\times X_p\right) \to G\times\varinjlim_p X_p
      \end{equation}
      (the colimits being taken in the category of topological spaces) is a homeomorphism. This can be seen by first realizing the direct limit as a quotient (in the topological sense \cite[\S 2]{mnk}) of a disjoint union, and then recalling that the product of a quotient map with (the identity morphism on) a compact Hausdorff space is again a quotient map \cite[Chapter 3, Exercise 11]{mnk}.

      The same line of reasoning applied to the left-hand factor $G$ in \Cref{eq:xtimesvarinjlim} then provides a homeomorphism
      \begin{equation}\label{eq:kpcolim}
        G\times X\cong \varinjlim_{K,p}K\times X_p,
      \end{equation}
      for $p$ ranging over the $C^*$ seminorms on $A$ and $K$ ranging over the compact subsets of $G$ (ordered by inclusion). Now, the universality property of the topological colimit \Cref{eq:kpcolim} makes the continuity of $G\times X\to X$ equivalent to that of each restriction
      \begin{equation*}
        K\times X_p\to X,
      \end{equation*}
      which we have seen already factors through some $X_q$. All in all, this reduces the continuity assertion to showing that every (co)restriction \Cref{eq:kpq} is continuous. To that end, consider nets
      \begin{equation*}
        g_{\lambda}\xrightarrow[\quad\lambda\quad]{} g\in K\quad\text{and}\quad \chi_{\lambda}\xrightarrow[\quad\lambda\quad]{} \chi\in X_p
      \end{equation*}
      The topology on each $X_q$ is nothing but the point-standard topology on $X_q$ regarded as a set of maps $A\to A_q\to \bC$, so the goal is to argue that for every $a\in A$ we have
      \begin{equation*}
        \chi_{\lambda}(g_{\lambda}^{-1}a)
        \xrightarrow[\quad\lambda\quad]{}
        \chi(g^{-1}a).
      \end{equation*}
      We may as well assume $q=\sup_K{}^g p$, as in the proof of part \Cref{item:kpq}, so that $g_\lambda$ and $g$, being in $K$, preserve that seminorm. But then, writing
      \begin{equation*}
        \chi_{\lambda}(g_{\lambda}^{-1}a) = \chi_{\lambda}(g_{\lambda}^{-1}a-g^{-1}a) + \chi_{\lambda}(g^{-1}a),
      \end{equation*}
      note that
      \begin{itemize}
      \item the first term is small for large $\lambda$, because the argument is small with respect to $q$ and $\chi_{\lambda}\in X_p$ are (norm-1) characters on $A_p$ (hence also $A_q$);
      \item while the second term is small for large $\lambda$ because $\chi_{\lambda}\to \chi$ in the point-standard topology. 
      \end{itemize}
    \end{enumerate}    
  \end{enumerate}
  This concludes the proof as a whole. 
\end{proof}

%%%%%%%%%%%%%%%%%%%%%%%%%%%%%%%%%%%%%%%%%%%%%%%%%%%%%%%%%%%%%%%%%%%%%%%%%%%%%
%%%%%%%%%%%%%%%%%%%%%%%%%%%%%%%%%%%%%%%%%%%%%%%%%%%%%%%%%%%%%%%%%%%%%%%%%%%%%

\section{Non-Commutative classifying space}\label{se:nus}

We refer the reader to~\cite[\S 14.1, 14.3, 14.4]{t-td08} for a detailed study of numerable principal bundles in topology. Recall that a numerable principal $G$-bundle $E\to E/G=:B$ is called {\em universal} if for any other numerable principal $G$-bundle $X\to X/G$ there exists a $G$-map $X\to E$ which is unique up to $G$-homotopy. All universal principal bundles are equivalent up to $G$-homotopy and one can prove that the space $B$ classifies all numerable principal $G$-bundles (see e.g.~\cite[Theorem~14.4.1]{t-td08}). The spaces $E$ and $B$ are called the~{\em universal $G$-space} and the~{\em classifying $G$-space} respectively. A~first model of a universal $G$-space for an arbitrary topological group $G$ was given by Milnor in~\cite[Section~2]{j-m56} as the topological join of infinitely many copies of $G$ with an appropriate topology (which is in general different than the quotient topology).  The choice of the topology is such that the diagonal action of $G$ on the infinite join is continuous and gives rise to a numerable principal $G$-bundle. One can show that the infinite join is contractible~(see e.g.~\cite[Proposition~14.4.6]{t-td08}). In fact, a numerable principal $G$-bundle is universal if and only if its total space is contractible~(\cite[Theorem~14.4.12]{t-td08}). All locally compact Hausdorff topological groups in this section are endowed with their natural quasi-topology.

%%%%%%%%%%%%%%%%%%%%%%%%%%%%%%%%%%%%%%%%%%%%%%%%%%%%%%%%%%%%%%%%%%%%%%%%%%%%%

\subsection{Cones of quasi-topological spaces}\label{subse:cone}

Let $X$ be a quasi-topological space. We define the {\em unreduced cone} of $X$ to be the quotient quasi-topological space
\begin{equation*}
\con(X):=(X\times[0,1])/(X\times \{0\})
\end{equation*}
(see~\cite[p.~5]{e-s63} for the explanation of the notation). We denote the quotient map by $\pi:X\times[0,1]\to\con(X)$.
If $G$ is a locally compact Hausdorff group acting quasi-continuously on $X$, then we define a $G$-action on $\con(X)$ by
\begin{equation}\label{actcon}
g[(x,t)]:=[(gx,t)],\qquad g\in G,\quad t\in[0,1],\quad x\in X.
\end{equation}
\begin{proposition}\label{prop:conqact}
The $G$-action on $\con(X)$ is quasi-continuous.
\end{proposition}
\begin{proof}
Let $\widetilde{\gamma}:G\times\con(X)\to\con(X)$ denote the action given by~\Cref{actcon}. From~\Cref{prodquot} we infer that the product quasi-topology on $G\times\con(X)$ is equivalent to the quotient quasi-topology induced by the map 
\begin{equation*}
{\rm id}\times \pi:G\times X\times [0,1]\longrightarrow G\times\con(X).
\end{equation*}
Hence, by~\Cref{quotlem}, to prove the claim it suffices to show that $\widetilde{\gamma}\circ({\rm id}\times \pi)$ is quasi-continuous. Observe that
  \begin{equation*}
    \begin{tikzpicture}[auto,baseline=(current  bounding  box.center)]
      \path[anchor=base] 
      (0,0) node (l) {$G\times X\times[0,1]$}
      +(3,.5) node (u) {$G\times \con(X)$}
      +(3,-.5) node (d) {$X\times[0,1]$}
      +(6,0) node (r) {$\con(X)$,}
      ;
      \draw[->] (l) to[bend left=6] node[pos=.5,auto] {$\scriptstyle {\rm id}\times\pi$} (u);
      \draw[->] (u) to[bend left=6] node[pos=.5,auto] {$\scriptstyle \widetilde{\gamma}$} (r);
      \draw[->] (l) to[bend right=6] node[pos=.5,auto,swap] {$\scriptstyle \gamma\times\id$} (d);
      \draw[->] (d) to[bend right=6] node[pos=.5,auto,swap] {$\scriptstyle \pi$} (r);
    \end{tikzpicture}
  \end{equation*}
where $\gamma:G\times X\to X$ is the quasi-continuous $G$-action on $X$. Since the bottom maps are quasi-continuous, we are done.
\end{proof}

\begin{proposition}
The {\em coordinate} maps
\begin{align}
{\rm t}:\con(X)\longrightarrow[0,1]:\quad &[(x,t)]\longmapsto t, \label{eq:tcor}\\
{\rm x}:{\rm t}^{-1}((0,1])\longrightarrow[0,1]:\quad &[(x,t)]\longmapsto x, \label{eq:xcor}
\end{align}
are quasi-continuous.
\end{proposition}
\begin{proof}
Let $\pi_X:X\times[0,1]\to X$ and $\pi_{[0,1]}:X\times[0,1]\to [0,1]$ denote the canonical surjections.

(1) Quasi-continuity of ${\rm t}$. We have to prove that ${\rm t}\circ\varphi$ is continuous for all $\varphi\in Q(K,\con(X))$ and all compact Hausdorff $K$. From the definition of quotient quasi-topology there exist a compact Hausdorff $K'$, a continuous surjection $\pi':K'\to K$, and $\varphi'\in Q(K',X\times[0,1])$ such that $\pi\circ\varphi'=\varphi\circ\pi'$. Since ${\rm t}\circ\varphi\circ\pi'=\pi_{[0,1]}\circ\varphi'$ is continuous by the definition of the product quasi-topology and from the fact that  $\pi'$ is a~quotient map, from~\Cref{quotlem} we infer that ${\rm t}\circ\varphi$ is continuous.

(2) Quasi-continuity of ${\rm x}$. Let $\iota:{\rm t}^{-1}((0,1])\to\con(X)$ denote the subspace inclusion map. We have to prove that ${\rm x}\circ\varphi\in Q(K,X)$ for all $\varphi\in Q(K,{\rm t}^{-1}((0,1]))$ and all compact Hausdorff $K$. Denoting $\psi:=\iota\circ\varphi$, there exist a compact Hausdorff $K'$, continuous surjection $\pi':K'\to K$, and $\psi'\in Q(K',X\times[0,1])$ such that $\pi\circ\psi'=\psi\circ\pi'$. Since ${\rm x}\circ\varphi\circ\pi'=\pi_X\circ\psi'$ is quasi-continuous, it follows from~\Cref{quotlem} that ${\rm x}\circ\varphi$ is continuous.
\end{proof}
Using the above maps, we observe that:
\begin{proposition}
If $X\in\cQ_{cH}$ then also $\con(X)\in\cQ_{cH}$.
\end{proposition}

We warn the reader that if $X$ is a topological space the quotient quasi-topology on $\con(X)$ need not be equivalent to the quasi-topology induced by the quotient topology.
\begin{example}
Let $X$ be any non-compact, countably compact, and separable Hausdorff space, e.g. the {\em Novak space}~\cite[Example~112]{countertop}, and consider $\con(X)$ with its quotient topology. Then by~\cite[Theorem~5]{gb-cone} there is a compact Hausdorff subset $K$ of $\con(X)$ such that there is no compact subset $L\subseteq X\times[0,1]$ such that $K\subseteq \pi(L)$. Observe that the inclusion $\iota:K\to\con(X)$ is in $C(K,\con(X))$ but it cannot be in $Q(K,\con(X))$, where we view $X$ as a quasi-topological space and $Q$ is the quotient quasi-topology induced by $\pi$.
\end{example}

\begin{proposition}\label{prop:contop}
If $X\in\cQ_{cH}$, then the topology on $\con(X)$ given by the equivalence of categories in \Cref{pr:chqt} agrees with the initial topology given by the maps~\Cref{eq:tcor} and \Cref{eq:xcor}.
\end{proposition}
\begin{proof}
Recall from~\cite[proof~of~Proposition~2.6]{phil1} that the considered topology on $\con(X)$ is defined as follows: $U\subseteq \con(X)$ is open if for every compact Hausdorff space $K$ and every $\varphi\in Q(K,\con(X))$, the set $\varphi^{-1}(U)$ is open in $K$. The distinguished family $\mathscr{F}_{\con(X)}$ of compacts subsets on $\con(X)$ is the set of images of all $\varphi\in Q(K,\con(X))$. By the definition of the quotient quasi-topology, we obtain that
\begin{equation*}
\mathscr{F}_{\con(X)}=\{L\subseteq\con(X)~|~\text{$L$ is compact and $L\subseteq\pi(L')$ for some compact $L'\subseteq X\times[0,1]$}\},
\end{equation*}
where on $X$ we also consider the topology given by~\Cref{pr:chqt}. 
Since, arguing as in the proof of~\cite[Lemma~3.1]{mt22} (here we use the assumption that $X$ is completely Hausdorff), one can show that the initial topology given by the maps~\Cref{eq:tcor} and~\Cref{eq:xcor} is generated by the family~$\mathscr{F}_{\con(X)}$, the conclusion follows.
\end{proof}

%%%%%%%%%%%%%%%%%%%%%%%%%%%%%%%%%%%%%%%%%%%%%%%%%%%%%%%%%%%%%%%%%%%%%%%%%%%%%

\subsection{Universal quasi-topological space}\label{subse:eg}

We consider a slight modification of the Milnor construction~\cite[Section~2]{j-m56}. 

\begin{definition}\label{def:eg}
  Let $G$ be a locally compact Hausdorff group. We define the quasi-space
  \begin{equation*}
    EG:=\left\{\left[(t_j,g_j)\right]_{j\in\mathbb{N}}\in\con(G)^\bN~\Bigg{|}~\sum_{j=0}^\infty t_j=1\right\},
  \end{equation*}
  equipped with the subspace quasi-topology.
\end{definition}

In the original Milnor construction, one assumes that only finitely many $t_j$ are not zero, so that the sum $\sum_jt_j$ is always finite. The difference between the two constructions is analogous to the difference between generalized partitions of unity and locally finite ones. The {\em coordinate} maps
\begin{equation*}
{\rm t}_i:EG\longrightarrow [0,1]:[(t_j,g_j)]_j\longmapsto t_i,\qquad {\rm g}_i:{\rm t}_i^{-1}((0,1])\longrightarrow G:[(t_j,g_j)]_j\longmapsto g_i\,,\qquad i\in\mathbb{N},
\end{equation*}
are quasi-continuous. Indeed, let ${\rm t}$ and ${\rm g}$ denote the maps~\Cref{eq:tcor} and~\Cref{eq:xcor} respectively and let 
$\pi_i:\con(G)^\bN\to\con(G)$ denote the projection on the $i$th component. Then ${\rm t}_i={\rm t}\circ\pi_i$ and ${\rm g}_i={\rm g}\circ\pi_i|_{{\rm t}_i^{-1}((0,1])}$. One can use the above maps to conclude that:
\begin{proposition}
If $G$ is a a locally compact Hausdorff group, then $EG\in\cQ_{cH}$.
\end{proposition}

The next theorem shows that our construction yields a universal numerable principal bundle (see the beginning of this section for a definition).

\begin{theorem}\label{th:eguni}
Let $G$ be a locally compact Hausdorff group. If on $EG$ we consider the topology given by~\Cref{pr:chqt}, then the action of $G$ on $EG$ gives rise to a universal numerable principal $G$-bundle.
\end{theorem}
\begin{proof}
  We will retrace the usual argument for universality, following the blueprint of \cite[discussion post Theorem 14.4.2]{t-td08} with the caveat that all maps in sight are $\cQ_{cH}$-morphisms.

  \begin{enumerate}[(I),wide]
  \item \textbf{: $EG$ is numerable.} Note that the continuous functions $\{{\rm t}_j\}_{j\in\mathbb{N}}$, define a~generalized partition of unity on $EG$. From~\cite[Lemma~13.1.7]{t-td08}, we know that one can construct a~locally finite partition of unity $\{\sigma_j\}_{j\in\mathbb{N}}$ such that ${\rm supp}\,\sigma_j\subseteq {\rm t}_j^{-1}((0,1])$ for all $j\in\mathbb{N}$. The covering $\{U_j:={\rm t}_j^{-1}((0,1])\}_{j\in\mathbb{N}}$ is trivializing, since we have the $G$-maps ${\rm g}_j:U_j\to G$.  
    
  \item \textbf{: Uniqueness of equivariant maps $\bullet\to EG$.} The claim, more precisely, is that for any numerable principal $G$-bundle $E\to B$ in $\cQ_{cH}$ any two $G$-equivariant maps $E\xrightarrow{f,g} EG$ are equivariantly homotopic. This follows from a quasi-topological adaptation of \cite[Proposition 14.4.4]{t-td08}, to whose proof there is very little to add; we remind the reader how the argument goes.

    Recalling our notation in \Cref{def:eg} for the coordinate description of $EG$, expand $f$ and $g$ as
    \begin{equation*}
      E\ni x
      \quad
      \xmapsto{\quad}
      \quad
      \left[(s_n(x),f_n(x))\right]_n
      \quad\text{and}\quad
      \left[(t_n(x),g_n(x))\right]_n
      \in EG
    \end{equation*}
    respectively. Assume next that we have homotopies
    \begin{equation}\label{eq:2hmtps}
      \begin{aligned}
        \left[(s_n,f_n)\right]_n
        \quad
        &\sim
          \quad
          \left[(s_1,f_1),\ 0,\ (s_2,f_2)\ ,\ 0\ \cdots\right]
          \quad\text{and}\\
        \left[(t_n,g_n)\right]_n
        \quad
        &\sim
          \quad
          \left[0,\ (t_1,g_1),\ 0,\ (t_2,g_2),\ 0\ \cdots\right]
      \end{aligned}      
    \end{equation}
    in $\cQ_{cH}$, i.e. realizations of the left- and right-hand sides as the restrictions to $\{0\}\times E$ and $\{1\}\times E$ respectively of $G$-equivariant $\cQ_{cH}$-morphisms $[0,1]\times E\to EG$. The map
    \begin{equation*}
      [0,1]\times E
      \ni
      (t,x)
      \xmapsto{\quad}
      \left[
        (ts_1,f_1),\ ((1-t)t_1,g_1),\ (ts_2,f_2)\ ,\ ((1-t)t_2,g_2)\ \cdots
      \right]
    \end{equation*}
    then effects a homotopy between the right-hand sides of \Cref{eq:2hmtps}, finishing the proof. It thus remains to argue that we indeed have homotopies \Cref{eq:2hmtps}, for which purpose it is sufficient, by the interchangeability of the two claims, to treat only the $f$ (i.e. top) branch.    

    Observe first that in general, for quasi-topological spaces $X,Y\in \cQ_{cH}$ chaining sequences
    \begin{equation*}
      \varphi_0
      \xrsquigarrow{\quad H_0\quad}
      \varphi_1
      \xrsquigarrow{\quad H_1\quad}
      \varphi_2
      \xrsquigarrow{\quad H_2\quad}
      \cdots
    \end{equation*}
    of homotopies between $\cQ_{cH}$-maps $X\xrightarrow{\varphi_n}Y$ functions precisely as usual: regarded as quasi-continuous maps
    \begin{equation*}
      \left[\frac {1}{n+1},\ \frac {1}{n}\right]\times X
      \xrightarrow{\quad H_n\quad}
      Y
    \end{equation*}
    agreeing appropriately, the $H_n$ aggregate into a map $(0,1]\times X\to Y$ by equipping
    \begin{equation*}
      (0,1]\times X
      =
      \left(\coprod_n \left[\frac 1{n+1},\frac 1n\right]\times X\right)
      \bigg/
      \left(\text{obvious equivalence relation}\right)
    \end{equation*}
    with the quotient quasi-topology mentioned in \Cref{exs:quasitops}\Cref{item:quot.qtop}. This observation in hand, we can now implement the top homotopy of \Cref{eq:2hmtps} by chaining together homotopies that remove progressively more 0s, starting with
    \begin{equation*}
      [0,1]\times E
      \ni
      (t,x)
      \xmapsto{\quad}
      \left[(s_1,f_1),\ (ts_2,f_2),\ ((1-t)s_2,f_2),\ (ts_3,f_3),\ ((1-t)s_3,f_3)\ \cdots\right]
      \in EG.
    \end{equation*}
    
  \item \textbf{: Existence of equivariant maps $\bullet\to EG$.} This time the claim is that numerable $G$-bundles $E\xrightarrow{\pi} B$ (again, in $\cQ_{cH}$) do in fact admit $G$-equivariant maps $E\to EG$. Here too the familiar argument (as presented in the proof of \cite[Proposition 14.4.5]{t-td08}, say) transports over to our quasi-topological framework virtually verbatim.
    \begin{itemize}[wide]
    \item On the one hand, given a $G$-equivariant \emph{countable} cover $\displaystyle E=\bigcup_{n\in \bZ_{\ge 0}}U_n$ trivializing $\pi$ and subordinating a $G$-invariant point-finite partition of unity $(\nu_n)_n$,
      \begin{equation*}
        E\ni x
        \xmapsto{\quad}
        \left[(\nu_n(x),f_n(x))\right]_n
      \end{equation*}
      for $G$-equivariant $U_n\xrightarrow{f_n}G$ witnessing the triviality of the $G$-spaces $U_n$. 

    \item On the other hand, \cite[Lemma 13.1.8 and Corollary 13.1.9]{t-td08} go through to leverage arbitrary partitions of unity into countable ones. 
    \end{itemize}
  \end{enumerate}
\end{proof}

%%%%%%%%%%%%%%%%%%%%%%%%%%%%%%%%%%%%%%%%%%%%%%%%%%%%%%%%%%%%%%%%%%%%%%%%%%%%%%%%%%%%%%

\subsection{Pro-$C^*$-relations revisited}\label{subse:relat}

In this subsection we construct the non-commutative universal space using universal pro-$C^*$-algebras of generators and relations that were first studied by Phillips~\cite[Section~1.3]{phil2}. Here we follow a~more general approach due to Loring~\cite[Section~3]{ta-l10}.

Let $\mathscr{G}$ be a set. A~{\em null pro-$C^*$-relation} on $\mathscr{G}$ is the category $\mathscr{G}\downarrow\cC^*_{pro}$ consisting of pairs $(A,j)$, where $A$ is a pro-$C^*$-algebra and $j:\mathscr{G}\to A$ is a mapping. A~morphism from $(A,j)$ to $(B,k)$ is a continuous $*$-homomorphism $\varphi:A\to B$ such that $\varphi\circ j=k$. The set $\mathscr{G}$ is called the set of {\em generators}.

\begin{remark}\label{re:comma}
  The `$\downarrow$' notation indicates the {\it comma category} construction of \cite[\S II.6]{mcl}. In that setup, the category $\mathscr{G}\downarrow\cC^*_{pro}$ would be $\mathscr{G}\downarrow\cat{forget}$ where
  \begin{equation*}
    \cC^*_{pro}\xrightarrow{\quad\cat{forget}\quad}\cat{Set}
  \end{equation*}
  is the forgetful functor. 
\end{remark}

Recall that if $A_\lambda$, $\lambda\in\Lambda$, is a family of pro-$C^*$-algebras, then $\prod A_\lambda$ is again a pro-$C^*$-algebra. Indeed, it is plainly a $*$-algebra, which can be endowed with the product pro-$C^*$-topology. We denote the canonical projections by $\xi_\lambda:\prod A_\lambda\to A_\lambda$. The following two definitions are due to Loring~\cite[Definition~3.11, Definition~3.12]{ta-l10}.

\begin{definition}\label{def:prorel}
  Let $\mathscr{G}$ be a set. A {\em pro-$C^*$-relation} on $\mathscr{G}$ is a full subcategory $\mathscr{R}$ of $\mathscr{G}\downarrow\cC^*_{pro}$ such that:
  \begin{enumerate}[{(R}1{)}]
  \item\label{item:r1} The pair $(\{0\},\chi)\in\mathscr{R}$, where $\{0\}$ is the zero pro-$C^*$-algebra and $\chi$ is the unique map $\mathscr{G}\to\{0\}$.
  \item\label{item:r2} Let $\varphi: A\hookrightarrow B$ be the inclusion of a closed $*$-subalgebra $A$ of a pro-$C^*$-algebra $B$ and let $f:\mathscr{G}\to A$ be a function. If $(B,\varphi\circ f)\in\mathscr{R}$ then $(A,f)\in\mathscr{R}$.
  \item\label{item:r3} Let $\varphi: A\to B$ be a continuous $*$-homomorphism and let $f:\mathscr{G}\to A$ be a function. If $(A,f)\in\mathscr{R}$ then $(B,\varphi\circ f)\in\mathscr{R}$.
  \item\label{item:r4} Let $\Lambda$ be a non-empty set and, for each $\lambda\in\Lambda$, let $f_\lambda:\mathscr{G}\to A_\lambda$ be a function. If $(A_\lambda,f_\lambda)\in\mathscr{R}$ for every $\lambda\in\Lambda$ then $(\prod A_\lambda,\prod f_\lambda)\in\mathscr{R}$.
  \end{enumerate}
\end{definition}

\begin{definition}\label{def:univloring}
  Let $\mathscr{R}$ be a pro-$C^*$-relation on a set $\mathscr{G}$. An object $(U,u)\in\mathscr{R}$ is called {\em universal} if the following conditions are satisfied:
  \begin{enumerate}
  \item[(U1)] Given a pro-$C^*$-algebra $A$, if $\varphi:U\to A$ is a continuous $*$-homomorphism then $(A,\varphi\circ u)\in\mathscr{R}$.
  \item[(U2)] Given a pro-$C^*$-algebra $A$, if $(A,f)\in\mathscr{R}$ then there is a unique $*$-homomorphism $\varphi:U\to A$ so that $f=\varphi\circ u$.
  \end{enumerate}
\end{definition}
By~\cite[Theorem~3.13]{ta-l10}, a universal pro-$C^*$-algebra exists (and is unique up to isomorphism) for every pro-$C^*$-relation $\mathscr{R}$ on a set $\mathscr{G}$.

\begin{example}\label{ex:sumtn1}
  Consider the following set of generators
  \begin{equation*}
    \mathscr{G}:=\{1,\; t_n~:~n\in\mathbb{N}\}.
  \end{equation*}
  Next, let $(A,f)\in \mathscr{G}\downarrow\cC^*_{pro}$. Then $(A,f)\in\mathscr{R}$ if and only if
  \begin{enumerate}
  \item[(1)] $f(1)=f(1)^*=f(1)^2$,\qquad $f(t_n)f(1)=f(1)f(t_n)=f(t_n)$,\quad $\forall~n\in\mathbb{N}$,
  \item[(2)] $p(f(t_n))\leq 1$, \quad $\forall~p\in S(A),~n\in\mathbb{N}$,
  \item[(3)] $\sum_{n=0}^\infty f(t_n)=f(1)$ in the pro-$C^*$-topology.
  \end{enumerate}
  It is straightforward to prove that $\mathscr{R}$ is a pro-$C^*$-relation. For instance, to show that $(2)$ and $(3)$ satisfy \ref{item:r4}, one could use~\cite[Lemma~3.6]{ta-l10}, which states that the pro-$C^*$-topology on $\prod A_\lambda$ is generated by the $C^*$-seminorms of the form
  \begin{equation*}
    {\rm max}\{p_0\circ\xi_{\lambda_j},\ldots,p_m\circ\xi_{\lambda_m}\},\qquad m\in\mathbb{N},
  \end{equation*}
  where $p_j\in S(A_{\lambda_j})$, $j=0,\ldots,m$. 
  One can prove that there is no universal $C^*$-algebra for the above relation. %above relation is not weakly admissible~(\cite[Definition~1.3.4]{phil2}).
  
  We will describe an object $(A,f)\in \mathscr{R}$ with {\it commutative} $A$. By \cite[Theorem~2.7]{phil1}, this entails specifying a completely Hausdorff quasi-space.
  The space is the ``infinite-dimensional simplex''
  \begin{equation*}
    \Delta_\infty:=\left\{(t_0,\ t_1,\ldots)\in[0,1]^{\bN}~\Big{|}~\sum_{i=0}^\infty t_i=1\right\}
  \end{equation*}
  embedded in the Hilbert cube $[0,1]^{\bN}$. The quasi-space structure is constructed virtually tautologically, so as to ensure the requisite universality property (always keeping in mind the duality of \Cref{th:antieq}). For compact Hausdorff $K$, the maps in $Q(K,\Delta_\infty)$ are precisely those $\varphi:K\to \Delta_\infty$ such that $t_i\circ\varphi$ is continuous for all $i\in\bN$, where $t_i$ denote the coordinate functions (this quasi-topology corresponds to the subspace product topology via~\Cref{pr:chqt}; note however that $\Delta_\infty$ is not closed in $[0,1]^\bN$). Since the sequence $(\sum_{i=0}^nt_i\circ\varphi)$ is clearly increasing and converges pointwise to $1$, by {\em Dini's} theorem~(see e.g.~\cite[Theorem~7.2.2]{dieu} for the metric case but observe that the same proof works for non-metric compact spaces), we obtain that
  \begin{equation*}
    \sum_{i=0}^{n} t_i\circ\varphi \xrightarrow[~n\to\infty~]{} 1\text{ uniformly on }K.
  \end{equation*}

Having equipped $A:=C(\Delta_{\infty})$ with its pro-$C^*$-algebra structure attached to the quasi-topology just described, the map
  \begin{equation*}
    \mathscr{G}\ni t_i\xmapsto{\quad f\quad}\left(\text{coordinate function }t_i\text{ on }\Delta_\infty\right)\in A
  \end{equation*}
  gives the desired pair $(A,f)\in \mathscr{R}$.
\end{example}

There is another, perhaps more easily guessable object $(A,f)\in \mathscr{R}$, also with commutative $A$, for the same relation $\mathscr{R}$ of \Cref{ex:sumtn1}:

\begin{example}\label{ex:uniondeltan}
  Consider the union $X$ of the chain
  \begin{equation*}
    \Delta_1\subset \Delta_2\subset \cdots\subset \Delta_n\subset\cdots
  \end{equation*}
  of finite-dimensional simplices
  \begin{equation*}
    \Delta_n = \left\{(t_0,\ \cdots,\ t_n)\in [0,1]^{n+1}\ |\ \sum t_i=1\right\},
  \end{equation*}
  with each inclusion adding an additional vertex. $X$ can then be equipped with the resulting colimit topology (so that it will be a {\it CW complex} \cite[preceding Example 0.1]{hatch_at}), and will produce a pro-$C^*$-algebra
  \begin{equation*}
    A:=C(X)\cong \varprojlim_n C(\Delta_n)
  \end{equation*}
  via \cite[Theorem~2.7]{phil1} again. Here too, there is a natural $f:\mathscr{G}\to A$, sending $t_n$ to the $n^{th}$ coordinate function in each $C(\Delta_m)$ (that coordinate function being 0 if $m<n$).

  The canonical map $X\to\Delta_{\infty}$ to the universal quasi-space of \Cref{ex:sumtn1} identifies $X$ with the (proper) union of the finite-dimensional simplices, so the two examples are genuinely different.
\end{example}

Next, let us rephrase \Cref{def:prorel} (that this is indeed nothing but a rephrasing, we will see later):

\begin{definition}\label{def:prorel-bis}
  Let $\mathscr{G}$ be a set and $\cC=\cC^*_{\bullet}$ one of the categories of (pro-)$C^*$-algebras in \Cref{not:cats} \Cref{item:allcast}.
  \begin{enumerate}[(1)]
  \item A {\it ($\mathscr{G}$-based) continuous $\cC$-relation} is a continuous subfunctor
    \begin{equation}\label{eq:ginpow}
      \Phi\lhook\joinrel\xrightarrow{\quad\iota\quad}(-)^{\mathscr{G}}
    \end{equation}
    of
    \begin{equation*}
      \cC\ni A\xmapsto{\quad} A^{\mathscr{G}}\in \cat{Set},
    \end{equation*}
    with $A^{\mathscr{G}}$ denoting the $\mathscr{G}$-indexed power in the category $\cC$, such that for every closed embedding $A\lhook\joinrel\xrightarrow{} B$ the square
    \begin{equation}\label{eq:embcar}
      \begin{tikzpicture}[auto,baseline=(current  bounding  box.center)]
        \path[anchor=base] 
        (0,0) node (l) {$\Phi(A)$}
        +(2,.5) node (u) {$\Phi(B)$}
        +(2,-.5) node (d) {$A^{\mathscr{G}}$}
        +(4,0) node (r) {$B^{\mathscr{G}}$}
        ;
        \draw[->] (l) to[bend left=6] node[pos=.5,auto] {$\scriptstyle $} (u);
        \draw[right hook->] (u) to[bend left=6] node[pos=.5,auto] {$\scriptstyle \iota_B$} (r);
        \draw[right hook->] (l) to[bend right=6] node[pos=.5,auto,swap] {$\scriptstyle \iota_A$} (d);
        \draw[->] (d) to[bend right=6] node[pos=.5,auto,swap] {$\scriptstyle $} (r);
      \end{tikzpicture}
    \end{equation}
    is a pullback. We refer to~\Cref{eq:embcar} as the {\it embedding-cartesian} condition (alternative language: \Cref{eq:ginpow} is embedding-cartesian).
  \item A relation as above is {\it finitely-continuous} if it preserves only {\it finite} limits.    
  \end{enumerate}
  Relations will typically be understood to be (fully) continuous unless otherwise stated. %explicitly flagged as finitely so.
\end{definition}
Note that~\Cref{eq:embcar} says something to the effect that the embedding \Cref{eq:ginpow} is a {\it cartesian natural transformation} \cite[Definition 4.1.1 c.]{leins_higher} when restricted to the subcategory of $\cC$ consisting of closed embeddings (rather than arbitrary morphisms). Continuity means, as usual \cite[\S V.4]{mcl}, that the functor preserves all (small) limits. The condition is pertinent:
\begin{itemize}
\item $\cC^*_{pro,\bullet}$ are complete even in the $\cQ$-enriched sense (see~\cite[Theorem~3.6]{catm-25});
\item while $\cC^*$, $\cC^*_1$, etc. are easily seen to also be complete: the equalizers are as in $\cC^*_{pro}$, while the products are the {\it bounded} tuples in the usual Cartesian product (and products and equalizers suffice \cite[Theorem 12.3 and Proposition 13.4]{ahs}).
\end{itemize}

To remove all potential ambiguity in the sequel (and justify the terminology coincidence), we first explain how \Cref{def:prorel-bis} is \Cref{def:prorel} re-tailored. 

\begin{proposition}\label{pr:samerel}
  Let $\mathscr{G}$ be a set.
  \begin{enumerate}[(1)]

  \item\label{item:pro-samerel} A $\cC^*_{pro}$-relation $\mathscr{R}$ in the sense of \Cref{def:prorel} gives rise to one in the sense of \Cref{def:prorel-bis} by setting
    \begin{equation}\label{eq:r2g}
      \Phi(A):=\{f:\mathscr{G}\to A\ |\ (A,f)\in \mathscr{R}\}.
    \end{equation}

    Conversely, a subfunctor \Cref{eq:ginpow} meeting the requirements of \Cref{def:prorel-bis} gives a $\cC^*_{pro}$-relation $\mathscr{R}$ in the sense of \Cref{def:prorel} by setting
    \begin{equation*}
      \mathscr{R}:=\{(f,A)\in \mathscr{G}\downarrow\cC^*_{pro}\ |\ f\in \Phi(A)\}.
    \end{equation*}    
    These two procedures are mutually inverse. 

  \item Similarly, a finitely-continuous $\cC^*$-relation in the sense of \Cref{def:prorel-bis} is equivalent to a $\cC^*$-relation in the sense of \cite[Definition 2.2]{ta-l10}. Under that correspondence, the (fully) continuous $\cC^*$-relations of \Cref{def:prorel-bis} are precisely the {\it compact} $\cC^*$-relations of \cite[Definition 2.3]{ta-l10}.
  \end{enumerate}
\end{proposition}
\begin{proof}
  This is all fairly routine (and an instance of the familiar correspondence between functors from a category to $\cat{Set}$ and {\it fibrations} over that category \cite[Theorem 8.3.1]{borc_hndbk-2}), so we only mention a~few key signposts, focusing on \Cref{item:pro-samerel}).

  To translate conditions \ref{item:r1} up to \ref{item:r4} into those of \Cref{def:prorel-bis}, for instance, note:
  \begin{itemize}
  \item \ref{item:r3} is just the functoriality of $\Phi$.
  \item \ref{item:r1} simply says that the functor $G:\cC^*_{pro}\xrightarrow{}\cat{Set}$ sends the terminal object $\{0\}$ to the terminal object $\{*\}$, i.e. preserves {\it empty } products.
  \item Similarly, \ref{item:r4} says that $\Phi$ preserves all other, non-empty products.
  \item As for \ref{item:r2}, it simultaneously ensures equalizer preservation and the embedding-cartesian condition of \Cref{def:prorel-bis}. For the former, for instance, consider an equalizer
    \begin{equation*}
      \begin{tikzpicture}[auto,baseline=(current  bounding  box.center)]
        \path[anchor=base] 
        (0,0) node (l) {$A$}
        +(-2,0) node (ll) {$A'$}
        +(2,0) node (r) {$B$}
        ;
        \draw[right hook->] (ll) to[bend left=0] node[pos=.5,auto] {$\scriptstyle j$} (l);
        \draw[->] (l) to[bend left=6] node[pos=.5,auto] {$\scriptstyle g$} (r);
        \draw[->] (l) to[bend right=6] node[pos=.5,auto,swap] {$\scriptstyle h$} (r);
      \end{tikzpicture}
    \end{equation*}
    in $\cC^*_{pro}$. For $\Phi$ constructed by \Cref{eq:r2g} we of course have $\Phi(gj)=\Phi(hj)$, $\Phi(j)$ is an embedding, and the fact that $(\Phi j)(\Phi A')$ {\it exhausts}
    \begin{equation*}
      \{x\in \Phi A\ |\ (\Phi g)(x)=(\Phi h)(x)\}
    \end{equation*}
    follows from \ref{item:r2}. 
  \end{itemize}
\end{proof}

With \Cref{def:prorel-bis} in place (and \Cref{pr:samerel} to connect it back to \cite[\S 3]{ta-l10}), the universal objects of \Cref{def:univloring} have the following interpretation:

\begin{lemma}\label{le:univrep}
  Let $\cC$ be one of the categories of \Cref{def:prorel-bis}. 
  
  Under the correspondences of \Cref{pr:samerel}, a universal object $(U,u)$ in the sense of \cite[Definitions 2.9 and 3.12]{ta-l10} is a realization of $\Phi:\cC\to \cat{Set}$ as a representable functor \cite[Definition 6.9]{ahs}: an object $U\in \cC$ together with a natural isomorphism
  \begin{equation}\label{eq:u2phi}
    \cC(U,-) \cong \Phi
  \end{equation}
  of functors $\cC^*_{pro}\xrightarrow{}\cat{Set}$.  \qedhere
\end{lemma}

\begin{remark}\label{re:repres}
  Since an isomorphism \Cref{eq:u2phi} will transport the embedding $\Phi\lhook\joinrel\xrightarrow{\iota}(-)^{\mathscr{G}}$ to
  \begin{equation*}
    \cC(U,-) \lhook\joinrel\xrightarrow{\quad}(-)^{\mathscr{G}},
  \end{equation*}
  the identity morphism on $U$ also provides us with a $\mathscr{G}$-tuple of elements in $U$, or, equivalently, a set map $\mathscr{G}\xrightarrow{f} U$ (and conversely, this data is sufficient to recover \Cref{eq:u2phi}). We will then also say that the relation $(\Phi,\iota)$ is represented by $f$ (or by $(U,f)$).
\end{remark}

Given all of the above, the existence \cite[Theorems 2.10 and 3.13]{ta-l10} of universal objects can be recovered as a representability statements:

\begin{theorem}\label{th:isrepr}
  Let $\cC$ be one of the categories of \Cref{def:prorel-bis} and $\mathscr{G}$ a set.

  A finitely-continuous $\mathscr{G}$-based $\cC$-relation in the sense of \Cref{def:prorel-bis} is representable if and only if it is fully continuous. 
\end{theorem}
\begin{proof}
  Representable functors are of course continuous, so one implication is immediate. 
  
  Conversely, assume a functor $\Phi$ as in \Cref{def:prorel-bis} is continuous. Its representability then follows provided $\Phi$ also satisfies the {\it solution-set condition} of \cite[\S \S V.6, Theorem 3]{mcl}: there is a {\it set} (as opposed to a proper class \cite[\S I.6]{mcl}) of objects $A_i\in \cC$ so that for every $A\in \cC$ and every $x\in \Phi(A)$ there is some morphism
  \begin{equation*}
    f:A_i\to A\text{ in }\cC\text{ with }x=(\Phi f)(y)\text{ for some }y\in \Phi(A_i).
  \end{equation*}
  To see this, observe that the $A$-component
  \begin{equation*}
    \Phi(A)\lhook\joinrel\xrightarrow{\quad\iota_A\quad}A^{\mathscr{G}}
  \end{equation*}
  of the embedding \Cref{eq:ginpow} realizes $x\in \Phi(A)$ as a $\mathscr{G}$-tuple of elements of $A$. It follows from the embedding-cartesian condition of \Cref{def:prorel-bis} that $x$ is an element of $\Phi(B)\le \Phi(A)$ for some closed (pro-)$C^*$-subalgebra $B\le A$ generated (as such; i.e. topologically) by a subset of $A$ of cardinality $\le |\mathscr{G}|$. 

  Since pro-$C^*$-algebras are by definition Hausdorff, we have \cite[\S 2.4]{juh_10}
  \begin{equation*}
    |A|\le \exp\exp\max\left(\aleph_0,|\mathscr{G}|\right)
  \end{equation*}
  for every (pro-)$C^*$-algebra $A$ generated by at most $|\mathscr{G}|$ elements. There is thus (up to isomorphism) only a {\it set} of such pro-$C^*$-algebras, and we can take those for our $\{A_i\}$. 
\end{proof}

\Cref{th:isrepr} provides, more or less, a correspondence between relations and objects of the various categories $\cC^*_{\bullet}$ under consideration. There are also matching {\it abelianization} procedures on the two sides of this picture:

\begin{definition}\label{def:ab}
  \begin{enumerate}[(1)]
  \item\label{item:abalg} The {\it abelianization} $A_{c}$ or $A_{ab}$ of an object $A\in \cC^*_{\bullet}$ is the completion in the appropriate category $\cC^*_{\bullet}$ of the quotient of $A$ by its closed $*$-ideal generated by the elements
    \begin{equation*}
      [a,a']:=aa'-a'a,\quad a,b.
    \end{equation*}
    This gives a functor
    \begin{equation*}
      \cC^*_{\bullet}\xrightarrow{\quad (-)_c=(-)_{ab}\quad}\cC^*_{\bullet,c}
    \end{equation*}
    (which is nothing but the identity in the uninteresting case when $\bullet$ already contains `c').
    
  \item\label{item:abrel} The {\it abelianization} $(\Phi,\iota)_{c}$ or $(\Phi,\iota)_{ab}$ of a $\cC^*_{\bullet}$-relation as in \Cref{def:prorel-bis} is the $\cC^*_{\bullet}$-relation consisting of the embedding
    \begin{equation*}
      \Phi_c\lhook\joinrel\xrightarrow{\quad\iota_{c}\quad}(-)^{\mathscr{G}}
    \end{equation*}
    where
    \begin{equation*}
      \Phi_c(A):=\left\{\text{tuples in }\Phi(A)\lhook\joinrel\xrightarrow{\iota}A^{\mathscr{G}}\text{ consisting of commuting elements of }A\right\}
    \end{equation*}
    and $\iota_c$ is the composition
    \begin{equation*}
      \Phi_c\lhook\joinrel\xrightarrow{\text{obvious inclusion}}\Phi\lhook\joinrel\xrightarrow{\quad\iota\quad}(-)^{\mathscr{G}}.
    \end{equation*}    
  \end{enumerate}
\end{definition}

The relation between the above two notions is precisely as expected:

\begin{proposition}\label{pr:abab}
  Let $(\Phi,\iota)$ be a continuous $\cC^*_{\bullet}$ relation, represented by $\mathscr{G}\xrightarrow{f} U$.

  The abelianization
  \begin{equation*}
    (\Phi,\iota)_c = (\Phi_c,\iota_c)
  \end{equation*}
  of \Cref{def:ab} \Cref{item:abrel} is then represented by
  \begin{equation*}
    \mathscr{G}\xrightarrow{\quad f\quad}U\xrightarrow{\quad}U_c,
  \end{equation*}
  where the right-hand map is the canonical morphism from $U$ to its abelianization in the sense of \Cref{def:ab} \Cref{item:abalg}. 
\end{proposition}
\begin{proof}
  This is a straightforward consequence of what it means for an object to represent a functor: a $\cC^*_{\bullet}$-morphism $U\to A$ corresponds to an element of
  \begin{equation*}
    \Phi_c\lhook\joinrel\xrightarrow{\quad}\Phi\lhook\joinrel\xrightarrow{\quad} (-)^{\mathscr{G}}
  \end{equation*}
  if and only if it sends the distinguished $\mathscr{G}$-tuple $(u_g)_{g\in \mathscr{G}}$ of $U$ into a commuting tuple of $A$. Because the $u_g$ are easily seen to generate $U$ as a $\cC^*_{\bullet}$-object, this is in turn equivalent to $U\to A$ factoring through a map $U_c\to A$. 
\end{proof}

%%%%%%%%%%%%%%%%%%%%%%%%%%%%%%%%%%%%%%%%%%%%%%%%%%%%%%%%%%%%%%%%%%%%%%%%%%%%%%%%%%%%%%

\subsection{Non-Commutative generalizations of $EG$ and $BG$}\label{subse:nceg}

Let us now describe the pro-$C^*$-algebra of quasi-continuous functions on the non-commutative space $E^+G$, which generalizes the quasi-space $EG$ from~\Cref{subse:eg}.
\begin{definition}\label{def:ce+g}
For a locally compact Hausdorff group $G$, $C(E^+G)$ is the object of $\cC^*_{pro,1}$ representing the continuous $\cC^*_{pro,1}$-relation $(\Phi,\iota)$ based on $\bZ_{\ge 0}\times C(\con(G))$ given by
  \begin{enumerate}[{(E}1{)}]
  \item\label{e1} For $A\in \cC^*_{pro,1}$, $\Phi(A)$ consists of morphisms
    \begin{equation*}
      \varphi_n:C(\con(G))\to A\text{ in }\cC^*_{pro,1},\quad n\in \bZ_{\ge 0};
    \end{equation*}
 so that
    \begin{equation*}
      \sum_{n=0}^\infty\varphi_n({\rm t})=1,\quad C(\con(G))\ni {\rm t}:=\text{height function}.
    \end{equation*}
  \item The inclusion
    \begin{equation*}
      \Phi(A)\lhook\joinrel\xrightarrow{\quad\iota\quad}A^{\bZ_{\ge 0}\times C(\con(G))}
    \end{equation*}
    is the obvious one in the set of arbitrary maps
    \begin{equation*}
      \bZ_{\ge 0}\times C(\con(G))\to A.
    \end{equation*}
  \end{enumerate}    
\end{definition}
First, we note that $\Phi$ is a pro-$C^*$-relation on $\mathscr{G}$. Indeed, the non-staighforward condition \ref{e1} follows from a reasoning analogous to the one outlined in~\Cref{ex:sumtn1}, which describes $C(E^+G)$ in the case of a trivial group.

\begin{remark}\label{rem:ceg}
  The pro-$C^*$-algebra $C(EG)$ for the quasi-space $EG$ of \Cref{def:eg} is obtained as in \Cref{def:ce+g}, with the added requirement that the images of the $\varphi_n:C(\con(G))\to A$ commute.
\end{remark}

An immediate consequence of \Cref{pr:abab}, as expected, the universal space $EG$ of \Cref{def:eg} is the ``classical shadow'' of $E^+G$:

\begin{corollary}[of~\Cref{pr:abab}]\label{cor:e+gab}
  For any $\cQ_{cH}$-group $G$ we have $C(EG)\cong C(E^+G)_{ab}$ as unital pro-$C^*$-algebras.
\end{corollary}
\begin{proof}
Per \Cref{def:ce+g} and \Cref{rem:ceg}, $C(EG)$ and $C(E^+G)$ are representing objects of relations, with the relation attached to the former being the abelianization (in the sense of \Cref{def:ab} \Cref{item:abrel}) of the relations corresponding to the former. \Cref{pr:abab} applies to yield the claim.
\end{proof}

Note, next, that the universality of $C(E^+G)$ affords us a natural $G$-action (no continuity claims for the moment):
\begin{itemize}
\item $G$ acts on $C(\con(G))$ by $\cC^*_{pro,1}$-automorphisms;
\item hence on the domains of the maps $\varphi_n:C(\con(G))\to A$ that make up $\Phi(A)$ in \Cref{def:ce+g};
\item and hence also on $\Phi(A)$ itself, since the $G$-action on $C(\con(G))$ fixes the height function ${\rm t}$;
\item so that $G$ operates on the functor $\Phi:\cC^*_{pro,1}\to \cat{Set}$ by natural automorphisms;
\item and hence also on the representing object $C(E^+G)$ by Yoneda \cite[Theorem 6.20]{ahs}. 
\end{itemize}

It is presumably not surprising that the expected continuity also holds:

\begin{proposition}\label{pr:gactse+g}
  The natural action of a $\cQ_{cH}$-group $G$ on $C(E^+G)$ is jointly continuous.
\end{proposition}
\begin{proof}
  The discussion preceding the statement already provides a monoid morphism
  \begin{equation}\label{eq:gactmaybediscont}
    G\to \cC^*_{pro,1}(C(E^+G),\ C(E^+G)),
  \end{equation}
  and we want to argue that it is quasi-continuous once the right-hand side is equipped with the following quasi-space structure: for compact Hausdorff $K$ the functions in $Q(K,C^*_{pro,1}(C(E^+G),C(E^+G))$ are those inducing continuous maps $K\times C(E^+G)\to C(E^+G)$.
  
  Now, $\cC^*_{pro,1}$ being $\cQ$-cotensored~\cite[Proposition~2.8]{catm-25}, there is a correspondence between quasi-continuous maps
  \begin{equation*}
    G\xrightarrow{\quad} \cC^*_{pro,1}(C(E^+G),\ C(E^+(G)))\in \cQ
  \end{equation*}
  and $\cC^*_{pro,1}$-morphisms
  \begin{equation}\label{eq:e+g2e+gg}
    C(E^+G)\xrightarrow{\quad} C(E^+G)^G
  \end{equation}
  (the latter being the pro-$C^*$-algebra of quasi-continuous maps $G\to C(E^+G)$). The defining universality property of $C(E^+G)$, in turn, identifies morphisms \Cref{eq:e+g2e+gg} with $\bZ_{\ge 0}$-tuples of morphisms $C(\con(G))\to C(E^+G)^G$ with the additional property that the images of the height functions add up to $1$. Now, the $G$-action \Cref{eq:gactmaybediscont} does induce maps
  \begin{equation*}
    C(\con(G))\xrightarrow{\quad f_n\quad}\cat{Set}(G,\ C(E^+G)),\quad n\in \bZ_{\ge 0}, 
  \end{equation*}
  and in order to conclude that these glue together to a morphism \Cref{eq:e+g2e+gg} it will be enough to argue that said maps take values in
  \begin{equation*}
    C(E^+G)^G\subset \cat{Set}(G,\ C(E^+G))
  \end{equation*}
  and are continuous. This, though, is immediate: $f_n$ factors as
  \begin{equation*}
    C(\con(G))\xrightarrow{\quad}C(\con(G))^G\subset \cat{Set}(G,\ C(\con(G)))\xrightarrow{\quad} \cat{Set}(G,\ C(E^+G)),
  \end{equation*}
  with the last arrow induced by the $n^{th}$ morphism $C(\con(G))\to C(E^+G)$ and the first map induced by the (jointly continuous) $G$-action on $C(\con(G))$. This implies both
  \begin{itemize}
  \item the claim that
    \begin{equation*}
      f_n C(\con(G))\subseteq C(E^+G)^G;
    \end{equation*}
  \item and the continuity of $f_n$, which has now been exhibited as the composition
    \begin{equation*}
      C(\con(G))
      \xrightarrow{\quad}
      C(\con(G))^G
      \xrightarrow{\quad}
      C(E^+G)^G
    \end{equation*}
    of two continuous maps. 
  \end{itemize}
  \end{proof}

Finally, let us define the analog of the classifying space in the non-commutative setting.
\begin{definition}\label{def:cb+g}
Let $G$ be a locally compact Hausdorff group and let $\alpha$ denote the quasi-continuous $G$-action on $C(E^+G)$ from~\Cref{pr:gactse+g}. The unital pro-$C^*$-algebra of the {\it non-commutative classifying space} $B^+G$ is defined as follows
\begin{equation*}
C(B^+G):=\{a\in C(E^+G)~|~\alpha_g(a)=a,~g\in G\}.
\end{equation*} 
\end{definition}

The notion of homotopy of $*$-homomorphisms of (unital) C*-algebras can be extended verbatim to (unital) pro-C*-algebras. Let us now generalize the {\it unital contractibility} of \cite[Definition~2.6]{dhn-20}).

\begin{definition}\label{def:homcontrpro}
  A unital pro-$C^*$-algebra $A$ is {\it unitally contractible} if the identity is homotopic in $\cC^*_{pro,1}$ to a map $A\to \bC\le A$.
\end{definition}

Properties of $EG$ familiar from the theory of principal bundles, such as the contractibility of $EG$ itself and that of the space of $G$-equivariant morphisms $X\to EG$ (for any $G$-space $X$) transport over to the non-commutative setup. The arguments are essentially those of \cite[Propositions 14.4.4 and 14.4.6]{t-td08}, slightly updated for the present context.

\begin{theorem}\label{th:gmaps2e+ghom}
  Let $G$ be a locally compact Hausdorff group and $\cC:=\cC^{*G}_{pro,1}$ the category of unital pro-$C^*$-algebras equipped with $G$-actions.

  Any two $G$-equivariant maps $C(E^+G)\to A$ are homotopic in $\cC$. 
\end{theorem}
\begin{proof}
  A morphism $f:C(E^+G)\to A$ as in the statement consists, by the defining universality property of the domain (\Cref{def:ce+g}), of a sequence of morphisms
  \begin{equation*}
    C(\con(G))\xrightarrow{\quad f_n\quad} A,\quad n\in \bZ_{\ge 0}
  \end{equation*}
  with
  \begin{equation*}
    \sum_n f_n({\rm t})=1;
  \end{equation*}
  the same goes for $g$ (with its own sequence $g_n$).

Let $A^X$ be the pro-$C^*$-algebra of quasi-continuous functions $X\to A$ equipped with the $G$-action induced by that of $A$. The desired homotopy is a morphism
  \begin{equation*}
    h:C(E^+G)\to A^I,\quad I:=[0,1]
  \end{equation*}
  specializing to $f$ and $g$ at the two endpoints $0,1\in I$. This, in turn, unpacks to a sequence
  \begin{equation*}
    h_n:C(\con(G))\to A^I,\quad n\in \bZ_{\ge 0}
  \end{equation*}
  of $G$-equivariant morphisms specializing to $f_n$ and $g_n$ respectively at $0$ and $1$ (keeping also in mind that $h_n({\rm t})$ must sum up to $1\in A^I$).

  The proof can now proceed as that of \cite[Proposition 14.4.4]{t-td08}. It will be convenient, given a~morphism $C(\con(G))\xrightarrow{\varphi} C$ and an $s\in I$, to write $s\varphi$ for the composition
  \begin{equation*}
    C(\con(G)) \xrightarrow{\quad} C(\con(G))\xrightarrow{\quad\varphi\quad}C
  \end{equation*}
  where the first arrow is dual to the scaling-by-$s$ map $\con(G)\to \con(G)$. In particular, for $s=0$ we simply write $0=0\varphi$. The first part of the proof of \cite[Proposition 14.4.4]{t-td08} then provides a~$G$-equivariant homotopy between
  \begin{itemize}
  \item $f=(f_n)$ and
    \begin{equation*}
      (f_0,\ 0,\ f_1,\ 0,\ \cdots);
    \end{equation*}
  \item and similarly, $g=(g_n)$ and
    \begin{equation*}
      (0,\ g_0,\ 0,\ g_1,\ \cdots).
    \end{equation*}
  \end{itemize}
  The second part of that same proof then gives a $t$-parametrized homotopy
  \begin{equation*}
    \left((1-t)f_0,\ tg_0,\ (1-t)f_1,\ tg_1,\ \cdots\right)
  \end{equation*}
  between these two latter morphisms. 
\end{proof}

And just as the argument employed in proving \cite[Proposition 14.4.4]{t-td08} also delivers the contractibility \cite[Proposition 14.4.6]{t-td08} of the usual classifying space $EG$, the proof of \Cref{th:gmaps2e+ghom} does double duty for the same purpose:

\begin{proposition}\label{pr:e+gcontractible}
  The unital pro-$C^*$-algebra $C(E^+G)$ is unitally contractible in the sense of \Cref{def:homcontrpro}.
\end{proposition}
\begin{proof}
  As in \Cref{th:gmaps2e+ghom}, one can first homotope the identity
  \begin{equation*}
    \id_{C(E^+G)} = (\iota_0,\ \iota_1,\ \cdots)
  \end{equation*}
  into
  \begin{equation*}
    (0,\ \iota_0,\ \iota_1,\ \cdots)
  \end{equation*}
  (even $G$-equivariantly), where
  \begin{equation*}
    C(\con(G))\xrightarrow{\quad\iota_n\quad} C(E^+G)
  \end{equation*}
  are the universal maps $C(E^+G)$ comes equipped with. The latter map is then homotopic via
  \begin{equation*}
    (t\cdot \mathrm{triv},\ (1-t)g_0,\ (1-t)g_1,\ \cdots),\quad t\in I
  \end{equation*}
  to $(1,\ 0,\ 0,\ \cdots)$, where
  \begin{equation*}
    C(\con(G))\xrightarrow{\quad \mathrm{triv}\quad}C(\con(G))
  \end{equation*}
  is the morphism dual to the map $\con(G)\xrightarrow{}\con(G)$ collapsing $\con(G)$ onto $1\in G$.
\end{proof}

%%%%%%%%%%%%%%%%%%%%%%%%%%%%%%%%%%%%%%%%%%%%%%%%%%%%%%%%%%%%%%%%%%%%%%%%%%%%%
%%%%%%%%%%%%%%%%%%%%%%%%%%%%%%%%%%%%%%%%%%%%%%%%%%%%%%%%%%%%%%%%%%%%%%%%%%%%%

\section{Non-Commutative principal bundles}\label{se:ncprin}

In this section, using the non-commutative universal space $E^+G$, we define $C^*$-algebraic principal bundles. First, we introduce the notion of a pro-$C^*$-algebraic principal bundle and then use the construction of the multipliers of the Pedersen ideal to adapt this definition to the category of $C^*$-algebras.

%%%%%%%%%%%%%%%%%%%%%%%%%%%%%%%%%%%%%%%%%%%%%%%%%%%%%%%%%%%%%%%%%%%%%%%%%%%%%

\subsection{Pro-$C^*$-algebraic principal bundles}

Given the construction~\Cref{def:ce+g}, we have a natural definition.
\begin{definition}\label{def:proprin}
A continuous action $\alpha:G\to{\rm Aut}(A)$ of a locally compact Hausdorff group on a~unital pro-$C^*$-algebra $A$, gives rise to a {\it pro-$C^*$-algebraic principal $G$-bundle} if there exists a~continuous unital $G$-equivariant $*$-homomorphism $C(E^+G)\to A$.
\end{definition}

We immediately observe that:
\begin{proposition}\label{pr:proprinab}
  A continuous action $\alpha:G\to {\rm Aut}(A)$  of a locally compact Hausdroff group $G$ on a unital commutative pro-$C^*$-algebra $A$ gives rise to a pro-$C^*$-algebraic principal bundle if and only if the induced $G$-action on ${\rm Spec}(A)$ (see~\Cref{th:antieq}) gives rise to a topological numerable principal bundle via the equivalence~\Cref{pr:chqt}.
\end{proposition}
\begin{proof}
  ($\Rightarrow$) Assume that $\alpha:G\to{\rm Aut}(A)$ defines a pro-$C^*$-algebraic principal bundle. Since $A$ is commutative, the $G$-equivariant morphism $C(E^+G)\to A$ in $\cC^*_{pro,1}$ factors through a $G$-equivariant morphism $C(E^+G)_{c}\to A$ in $\cC^*_{pro,c,1}$. From~\Cref{cor:e+gab} and~\Cref{th:antieq}, we obtain a $G$-equivariant morphism ${\rm Spec}(A)\to EG$ in $\cQ_{cH}$. Combining~\Cref{pr:chqt} and~\Cref{th:eguni}, we conclude that we have a $G$-equivariant continuous map ${\rm Spec}(A)\to EG$, which is equivalent with ${\rm Spec}(A)$ being a total space of a numerable principal $G$-bundle.

  ($\Leftarrow$) Assume that the induced $G$-action makes ${\rm Spec}(A)$ a total space of a numerable principal $G$-bundle. Like before, it is equivalent with the existence of a $G$-equivariant morphism ${\rm Spec}(A)\to EG$ in~$\cQ_{cH}$, which in turn is equivalent with a $G$-equivariant morphism $C(E^+G)_c\to A$ in $\cC^*_{pro,c,1}$. To obtain the claim, we compose this morphism with the canonical $G$-equivariant morphism $C(E^+G)\to C(E^+G)_c$ in $\cC^*_{pro,1}$.
\end{proof}

The result above ensures that~\Cref{def:proprin} recovers all classical numerable principal bundles for which the total spaces are objects in $\cat{Top}_{cH,\mathscr{F}}$.

\begin{example}
The $C^*$-algebra $C(S^n_+)$ of the free $n$-sphere $S^n_+$ (\cite[Definition~2.1]{bg10}) is the universal unital $C^*$-algebra generated by self-adjoint elements $x_i$, $i=0,1,\ldots,n$, subject to the relation
\begin{equation*}
\sum_{i=0}^n x_i^2=1.
\end{equation*}
There is a natural antipodal action $\alpha_n$ of $\bZ/2\bZ=\{1,-1\}$ on each $C(S^n_+)$ given on generators by 
\begin{equation*}
x\mapsto -x_i, \qquad i=0,1,\ldots,n.
\end{equation*} 
Consider the pro-$C^*$-algebra 
\begin{equation*}
C(S^\infty_+):=\varprojlim_nC(S^n_+)
\end{equation*}
with the inverse limit action $\alpha:=\varprojlim_n\alpha_n$ (see~\cite[Definition~5.1]{nc-p89}). Let $\gamma$ be the generator of the $C^*$-algebra $C(\bZ/2\bZ)$ and consider the following function $\sqrt{\rm t}\otimes\gamma:\con(\bZ/2\bZ)\to[0,1]:[(\pm 1,t)]:=\pm\sqrt{t}$. As in~\cite[Example~3.10]{ghtw-18}, the assignments
\begin{align*}
\varphi_0:C(\con(\bZ/2\bZ))\ni\sqrt{\rm t}\otimes \gamma&\mapsto (1,x_1,x_1,x_1\ldots)\in C(S^\infty_+),\\
\varphi_1:C(\con(\bZ/2\bZ))\ni\sqrt{\rm t}\otimes \gamma&\mapsto (0,x_2,x_2,x_2,\ldots)\in C(S^\infty_+),\\
\varphi_2:C(\con(\bZ/2\bZ))\ni\sqrt{\rm t}\otimes \gamma&\mapsto (0,0,x_3,x_3,\ldots)\in C(S^\infty_+),\\
&~\vdots
\end{align*}
define $\bZ/2\bZ$-equivariant unital $*$-homomorphisms $\varphi_n:C(\con(\bZ/2\bZ))\to C(S^\infty_+)$, $n\in\bZ_{\geq 0}$, such that
\begin{equation*}
\sum_{n=0}^\infty\varphi_n({\rm t})=1.
\end{equation*}
Hence, we get a $\bZ/2\bZ$-equivariant unital continuous $*$-homomorphism $C(E^+\bZ/2\bZ)\to C(S^\infty_+)$.
\end{example}

Every unital $C^*$-algebra is a unital pro-$C^*$-algebra. In this case, our definition can be compared with the notion of the local-triviality dimension~(compare with \cite{ghtw-18}):
\begin{definition}\label{def:dimlt}
Let $G$ be a compact Hausdorff group acting on a unital $C^*$-algebra $A$. We say that the {\it local-triviality dimension} of the $G$-$C^*$-algebra $A$ equals $n\in\bZ_{\geq 0}$, denoted $\dim_{\rm LT}^G(A)=n$, if $n$ is the smallest number such that there exist $G$-equivariant $*$-homomorphisms $\rho_n:C(\con(G))\to A$ such that $\sum_{i=0}^n\rho_n({\rm t})=1$. We write $\dim^G_{\rm LT}(A)=\infty$ if no such $n$ exists.
\end{definition}
It is immediate that:
\begin{proposition}\label{pr:dimprin}
If an action of a compact Hausdorff group $G$ on a unital $C^*$-algebra $A$ has finite local-triviality dimension, then it defines a pro-$C^*$-algebraic principal bundle.
\end{proposition}
The converse to \Cref{pr:dimprin} is, in part, the subject of future work. 

%%%%%%%%%%%%%%%%%%%%%%%%%%%%%%%%%%%%%%%%%%%%%%%%%%%%%%%%%%%%%%%%%%%%%%%%%%%%%

\subsection{$C^*$-algebraic principal bundles}\label{subse:multped}

In the present subsection, we deal with actions on non-unital $C^*$-algebras. First, we recall the construction of the algebra of multipliers on the Pedersen ideal. We define the following sets
\begin{equation*}
  K_0^+:=\{a\in A_+~:~ab=a\text{ for some $b\in A_+$}\},
\end{equation*}
\begin{equation*}
  K^+:=\left\{a\in A_+~:~\text{there exists $a_0\,,\ldots, a_n\in K_0^+$ such that }a\leq \sum_{i=0}^na_n\right\}.
\end{equation*}
In \cite[discussion preceding Theorem 1.3]{ped_meas-1}, Pedersen introduced a minimal dense ideal $K_A:={\rm span}\,K^+$ of $A$.  Let us denote its algebra of (not necessary bounded) multipliers by $\Gamma(K_A)$. It is a unital $*$-algebra with natural operations induced from $K_A$. Bounded multipliers of $K_A$ correspond to the multiplier algebra $M(A)$. In \cite[Theorem 4]{phil_multiplier}, Phillips proved the following isomorphism of unital $*$-algebras
\begin{equation*}
  \Gamma(K_A)\cong\varprojlim_{k\in(K_A)_+} M(I_k)\cong\varprojlim_\lambda M(I_{e_\lambda}).
\end{equation*}
Here $(K_A)_+$ denotes the positive elements in $K_A$, $I_k$ is the closed two-sided ideal generated by the element $kk^*$, $M(I_k)$ is its multiplier algebra, and $(e_\lambda)$ is an approximate unit of $A$ contained in $K_A$. The canonical $*$-homomorphisms 
\begin{equation*}
\pi_k:\Gamma(K_A)\to M(I_k):(S,T)\mapsto (S_k,T_k)
\end{equation*}
are given by restriction, namely $S_k:=S|_{I_k}$ and $T_k:=T|_{I_k}$. The above isomorphism defines a pro-$C^*$-topology on $\Gamma(K_A)$ via the inverse limit of norm topologies on the multiplier algebras $M(I_k)$. 

Observe that $\Gamma(K_A)$ is commutative if and only if $A$ is commutative. In the commutative case, i.e. when $A=C_0(X)$ for some locally compact Hausdorff space $X$, we have that $K_A$ is isomorphic with the $*$-algebra of compactly supported functions $C_c(X)$. Consequently, $\Gamma(K_A)$ is isomorphic with the $*$-algebra of all continuous functions $C(X)$.

%\begin{remark}\label{re:alltopsamecomm}
%  For commutative $A=C_0(X)$ the pro-$C^*$, $\kappa$ and {\it strict} (\cite[Theorem 7]{phil_multiplier}, \cite[Definition 3.13]{phil1}) topologies on $\Gamma(K_A)\cong C(X)$ all coincide with the topology of uniform convergence on compact sets (i.e. the compact-open topology). Indeed, the individual $A/I(X-C)$ in~\cite[Theorem~7]{phil_multiplier} are unital and isomorphic to $C(C)$ (continuous functions on compact subsets $C$ of the spectrum).
%\end{remark}

Assume now that $A$ carries a (jointly) continuous action of a locally compact Hausdroff group~$G$. It is clear that $\alpha$ can be extended to a $G$-action on $\Gamma(K_A)$. In general, this extended action will not be continuous, which prompts the following piece of notation
\begin{equation*}
\Gamma(K_A)_G:=\{a\in\Gamma(K_A)~|~(g,a)\mapsto \alpha_g(a)\text{ is (jointly) continuous in the pro-$C^*$-topology}\}.
\end{equation*}

\begin{definition}\label{def:cprin}
An action $\alpha:G\to {\rm Aut}(A)$ of a locally compact Hausdroff group $G$ on a $C^*$-algebra~$A$ gives rise to a~{\it $C^*$-algebraic principal $G$-bundle} if there exists a $G$-equivariant unital continuous $*$-homomorphism $C(E^+G)\to\Gamma(K_A)_G$.
\end{definition}
\begin{remark}
Observe that considering $\Gamma(K_A)_G$ is not restrictive. Indeed, the image of any $G$-equivariant continuous $*$-homomorphism $C(E^+G)\to\Gamma(K_A)$ is contained in $\Gamma(K_A)_G$.
\end{remark}

We conclude the paper with various classes of examples of $C^*$-algebraic principal bundles.
\begin{examples}
\begin{enumerate}[(1)]
\item As per~\Cref{pr:proprinab}, if $X$ is a locally compact Hausdorff $G$-space (meaning that $G$ is also locally compact Hausdorff), then $X\to X/G$ is a numerable principal $G$-bundle if and only if the induced action $\alpha:G\to {\rm Aut}(C_0(X))$ defines a $C^*$-principal $G$-bundle.
\item Let $G$ be compact Hausdorff and let $A$ be unital. Since then $\Gamma(K_A)\cong A$,~\Cref{def:cprin} coincides with~\Cref{def:proprin}. Hence, whenever $\dim_{\rm LT}^G(A)<\infty$ (see~\Cref{def:dimlt}), we get an example of a $C^*$-algebraic principal $G$-bundle.
\item Let $G$ be locally compact Hausdorff acting on a $C^*$-algebra $A$ and consider the $\kappa$-topology~\cite[\S 3.6]{lt_pedersen} on $\Gamma(K_A)$. Then we say that $A$ is a {\it locally trivial $G$-$C^*$-algebra}~\cite[Definition~4.12]{mt22} if there exist $G$-equivariant unital $\kappa$-continuous $*$-homomorphisms $\rho_n:C(\con(G))\to\Gamma(K_A)$ such that $\sum_{n=0}^\infty\rho_n({\rm t})=1$ in the $\kappa$-topology. Since the $\kappa$-topology is weaker than the pro-$C^*$-topology on $\Gamma(K_A)$, each locally trivial $G$-$C^*$-algebras produces a $C^*$-algebraic principal $G$-bundle.
\item Let $G$ be a second-countable torsion-free locally compact Hausdorff group. Following~\cite[Proposition~4.14]{mt22}, we note that in this case the {\it proper $G$-algebras} of Guenter, Higson, and Trout~\cite[Definition~1.3]{ght00} also define $C^*$-algebraic principal $G$-bundles.
\item Let $G$ be a locally compact Hausdorff abelian group acting locally unitarily on a $C^*$-algebra $A$ (see~\cite{rw88}). Again, following the discussion above~\cite[Remark~4.9]{mt22}, we note that dual actions on crossed products $G\ltimes A$ define $C^*$-algebraic principal $G$-bundles.
\end{enumerate}
\end{examples}

%%%%%%%%%%%%%%%%%%%%%%%%%%%%%%%%%%%%%%%%%%%%%%%%%%%%%%%%%%%%%%%%%%%%%%%%%%%%%%%%%%%%%%
%%%%%%%%%%%%%%%%%%%%%%%%%%%%%%%%%%%%%%%%%%%%%%%%%%%%%%%%%%%%%%%%%%%%%%%%%%%%%%%%%%%%%%

%\bibliographystyle{plain}
%\bibliography{universalrefs}{}

% BEGIN INSERTED BBL (uni_curr-xv2.bbl)
\def\polhk#1{\setbox0=\hbox{#1}{\ooalign{\hidewidth
  \lower1.5ex\hbox{`}\hidewidth\crcr\unhbox0}}}
  \def\polhk#1{\setbox0=\hbox{#1}{\ooalign{\hidewidth
  \lower1.5ex\hbox{`}\hidewidth\crcr\unhbox0}}}
  \def\polhk#1{\setbox0=\hbox{#1}{\ooalign{\hidewidth
  \lower1.5ex\hbox{`}\hidewidth\crcr\unhbox0}}}
  \def\polhk#1{\setbox0=\hbox{#1}{\ooalign{\hidewidth
  \lower1.5ex\hbox{`}\hidewidth\crcr\unhbox0}}}
  \def\polhk#1{\setbox0=\hbox{#1}{\ooalign{\hidewidth
  \lower1.5ex\hbox{`}\hidewidth\crcr\unhbox0}}}

% END INSERTED BBL

\Addresses

\end{document}